\documentstyle{article}
\begin{document}

\title{A Markov Property For Set-Indexed Processes}

\author{R.M. Balan\thanks{Corresponding author. 
Department of Mathematics and Statistics, University of Ottawa, P.O.
Box 450, Stn. A, Ottawa, ON, K1N 6N5, CANADA, Tel: (613)
260-0683. {\em E-mail address: s1412643@matrix.cc.uottawa.ca}}
\thanks{This paper is based on a portion of the first author's
doctoral thesis. Research supported by a scholarship from the Natural
Sciences and Engineering Research Council of Canada and by an Ontario
Graduate Scholarship in Science and Technology.}     
\and   
B.G. Ivanoff\thanks{Department of
Mathematics and Statistics, University of Ottawa, Ottawa, ON, CANADA.
{\em E-mail address: givanoff@science.uottawa.ca}} \thanks{Research
supported by a grant from the Natural Sciences and Engineering
Research Council of Canada.} }

\date{\today}

\maketitle

\begin{abstract}
We consider a type of Markov property for set-indexed processes which
is satisfied by all processes with independent increments and which
allows us to introduce a transition system theory leading to the
construction of the process. A set-indexed generator is defined such
that it completely characterizes the distribution of the process. 
\end{abstract}

\noindent {\em Keywords}: set-indexed process; Markov property;
transition system; generator.

\section{Introduction}

The Markov property is without doubt one of the most appealing
notions that exists in the classical theory of stochastic processes and
many processes modelling physical phenomena enjoy it. Attempts have
been made to generalize this concept to processes where the index
set is not totally ordered.

However, to define a Markov property for processes indexed by an
uncountable partially ordered set is not a straightforward task for
someone who has as declared goals to prove that: (i) all processes
with independent  increments possess this property; (ii) there
exists a systematic procedure which allows us to construct a general
process which enjoys this property; and (iii) we can define a
generator which completely characterizes the finite dimensional
distributions of the process. (The case of Markov processes indexed by
discrete partially ordered sets was considered by various authors
\cite{cairoli-dalang96}, \cite{lawler-vanderbei81} for solving
stochastic optimal control problems; as our framework is different,
we will not discuss this case here.) In the present paper we will
define a new type of Markov property for processes indexed by a semilattice of
sets, which will attain these three goals. 

This paper represents a successful manner of approaching a problem to 
which much effort has been dedicated in the literature to date. In fact it will be seen that our definition is completely analogous to the classical definition of the Markov property on ${\bf R}$, and our approach avoids many of the technical problems which arise with other definitions of Markov properties.
We begin now to describe the two levels of generality of this problem.

The first level deals with the case when the index set 
is the Euclidean space $[0,1]^{2}$ (or, more generally $[0,1]^{d}$)
and therefore it inherits 
the extra structure introduced by the total ordering of the coordinate 
axes. In this case, there are at least three types of Markov properties
which can be considered: the sharp Markov property, the germ Markov 
property, and the `simultaneously vertical and horizontal' Markov 
property.

The basic idea behind the sharp Markov property (first introduced in 
1948 by L\'{e}vy \cite{levy48}) was to consider as the 
history of the process at location $z$ in the plane or space, all the 
information that we have about the values of the process inside the 
rectangle $[0,z]$ ; all the information about the values of 
the process outside the rectangle was regarded as `future'; and the 
past and the future should be independent given the values of the
process on the boundary of the rectangle. In 1984, Russo \cite{russo84} 
proved that all processes with independent increments are sharp 
Markov with respect to all finite unions of rectangles.
The next step was to see if we can replace the rectangles with more 
general sets; in other words
a two-parameter process $(X_{z})_{z \in [0,1]^{2}}$ is said to have 
{\bf the sharp Markov property} with respect to a set $A$ if the $\sigma$-
fields ${\cal F}_{A}$ and ${\cal F}_{A^{c}}$ are conditionally 
independent given ${\cal F}_{\partial A}$ where ${\cal F}_{D}= 
\sigma(X_{z}; z \in D)$ for any set $D \subseteq [0,1]^{2}$. 
The advantage of this definition is that it does not rely on the
partial order of the space. Unfortunately, it turned out that processes 
having this property are difficult to handle. 
To attain goal (i) the following question had to be answered: what is 
the largest class of sets for which all processes with independent 
increments are sharp Markov? (To answer this question the Gaussian 
and the jump parts have been treated separately.) Goals (ii) and
(iii) do not seem to be even specified anywhere in the literature. 

(In 1976 Walsh \cite{walsh76} showed that the Brownian sheet fails 
to have the sharp Markov property
with respect to a very simple set, the triangle with
vertices $(0,0),(0,1)$ and $(1,0)$. This led to the conclusion that,
instead of the sharp $\sigma$-field ${\cal F}_{\partial A}$ one
has to consider a larger one, called the germ $\sigma$-field,
which is defined as ${\cal G}_{\partial A}=\cap_{G} {\cal
F}_{G}$ where the intersection is taken over all open sets $G$
containing $\partial A$. The new Markov property, for which
${\cal F}_{A}$ and ${\cal F}_{A^{c}}$ are conditionally
independent given ${\cal G}_{\partial A}$, was called 
{\bf the germ Markov property} and it was first introduced by McKean 
\cite{mckean63}. In 1980 Nualart \cite{nualart80} showed that the
Brownian sheet satisfies the germ Markov property with respect
to all open sets.)

The complete answer was given in 1992 in the form
of two thorough papers written by Dalang and Walsh: in
\cite{dalang-walsh92-br} it is shown that the Brownian sheet
satisfies the sharp Markov property with respect to all
domains whose boundaries are singular curves of
bounded variation; on the other hand, the main
result of \cite{dalang-walsh92-levy} states that for jump
proceses with independent increments which have only
positive jumps, the sharp Markov property holds with respect to all
domains. 

The third Markov property defined for two-parameter processes 
was introduced in 1979 by Nualart and Sanz 
\cite{nualart-sanz79} and it was first studied in the
Gaussian case. Korezlioglu, Lefort and Mazziotto \cite{KLM}
generalized it and proved that any process satisfying this property
is sharp Markov with respect to any finite union of
rectangles and germ Markov with respect to any relatively
convex domain. The basic idea in this third Markov property is the
separation of parameters, that is, the simultaneous definition
of a horizontal and vertical Markov property.
It is an easy exercise to show that any process with independent 
increments satisfies this property. The general construction of 
a process satisfying this property, which corresponds to a certain
transition semigroup, was made by Mazziotto \cite{mazziotto88}
(the trajectories of this process have also nice regularity 
properties).

The second level of generality deals with the case of processes
indexed by a collection ${\cal A}$ of closed connected subsets of a
compact  space. The collection ${\cal A}$ is partially ordered by
set-inclusion and has the additional structure of a semilattice,
being  closed under arbitrary intersections. This new approach in the
modern  theory of stochastic processes indexed by partially ordered
sets (which is known in the literature as the theory of set-indexed
processes) was initiated and developed by Ivanoff and Merzbach for
the martingale case \cite{ivanoff-merzbach00}. (Other authors
\cite{adler-feigin84},\cite{bass-pyke84} were 
interested in processes indexed by Borel subsets of $[0,1]^{d}$, but
the approach taken by Ivanoff and Merzbach introduces tools which
permit the development of a general  theory.) If we identify this
class with the class of the rectangles $[0,z],z \in [0,1]^{d}$ and we
denote by
$X_{z}$ the value of the set-indexed process at the rectangle
$[0,z]$, then $(X_{z})_{z}$ will be a $d$-dimensional process;
hence, we can view this theory as a generalization of the theory
of multiparameter processes. 

In \cite{ivanoff-merzbach98} both sharp Markov and (`simultaneously 
vertical and horizontal') Markov properties have been introduced
in the set-indexed framework, in such a manner 
that they generalize the corresponding notions defined in the
multiparameter case. In fact, in this new set-up, the Markov property 
implies the sharp Markov property and there are certain circumstances 
when they are actually equivalent. Moreover, a highly non-trivial 
result of \cite{ivanoff-merzbach98} shows that all processes with 
independent increments are Markov (in the specified sense), so that 
goal (i) is attained; however, the authors of \cite{ivanoff-merzbach98}
do not consider goals (ii) and (iii) and they do not attempt to
construct a general Markov process.

In this paper we will consider another type of Markov property
for set-indexed processes, which is less general than the Markov 
property introduced by Ivanoff and Merzbach, but which has the merit
of attaining the three goals (i), (ii) and (iii). 
This Markov property will be called {\bf the set-Markov property}
and we will show that it implies the sharp Markov property.
The definition requires that the value 
$X_{A \verb2\2 B}=X_{A \cup B}-X_{B}$ of
an additive set-indexed process $(X_{A})_{A \in {\cal A}}$ over the 
increment $A \verb2\2 B$
be conditionally independent of the history ${\cal F}_{B}:
= \sigma(\{X_{A'};A' \in {\cal A},A' \subseteq B \})$
given the present status $X_{B}$. We believe that this is a
natural definition because:
\begin{description}
\item[(a)] it captures the essence of the
Markov property in terms of the increments, allowing us to have
a perfect analogy with the classical case;
\item[(b)] all processes with independent increments are (trivially)
set-Markov; 
\item[(c)]  a  process  is set-Markov if and only if it
becomes Markov in the classical sense when it is transported by a
`flow'. (A `flow' is an increasing function mapping a bounded
interval of the real line into the collection of all finite unions
of  sets in ${\cal A}$.) This property provides us with the means to
define a generator for a set-Markov process. 
\end{description}

The paper is organized as follows:

In Section 2 we define the general framework which is used
in the theory of set-indexed processes. The main structure that we need is a semilattice 
of sets which leads us to a semialgebra and an algebra of sets.
All the processes that we consider are assumed to be a.s. finitely additive on this
algebra.

In Section 3 we define the set-Markov property in terms of the increments of the process.
We will prove that this property has two equivalent formulations. The first one (Proposition 1) requires that the `future' behaviour of the process at a set $B' (\supseteq B)$ be conditionally independent of its history ${\cal F}_{B}$ knowing its present status at the set $B$. The second formulation (Proposition 2) requires that the (one-parameter) processes $X^{f}$ obtained as `traces' of the original set-indexed process $X$ along the paths of all flows $f$ be Markov in the classical sense; this allows us to make use of the rich theory that exists for Markov processes indexed by the real line.

In Section 4 we turn to the question of existence of a set-Markov process and we prove that we can construct such a process if we specify the laws $Q_{BB'}$ that characterize the transition from a set $B$ to a set $B' (\supseteq B)$. A set-Markov process with a given transition system ${\cal Q}:=(Q_{BB'})_{B \subseteq B'}$ is called `${\cal Q}$-Markov' and can be constructed as soon as we specify its finite dimensional distributions (by Kolmogorov's existence theorem). However, the formula that we obtain for the finite dimensional distribution over some $k$-tuple $A_{1}, \ldots,A_{k}$ of sets turns out to be dependent on the ordering of 
the sets; hence we have to impose a natural consistency condition (Assumption 1) which makes the finite dimensional distribution `invariant' under permutations. An important consequence of this condition is that the finite dimensional distributions of a ${\cal Q}$-Markov process are also additive (Lemma 5).

In Section 5 we define the generator of
a ${\cal Q}$-Markov process $X$ as the class $\{{\cal G}^{f}; f \in {\cal S}\}$, where ${\cal G}^{f}$ is the generator of the (one-parameter) Markov process $X^{f}$ and ${\cal S}$ is a large enough (uncountable) collection of flows. In other words, at any set $B$ there are infinitely many generators depending on which `direction' we approach this set along the path of a flow $f$ in ${\cal S}$ (note that in the classical case, the class ${\cal S}$ of flows can be taken to have a single element). 

In this section we address the question of existence of a ${\cal Q}$-Markov process given its generator. More precisely, we start with a collection $\{ {\cal G}^{f}; f \in {\cal S}\}$, each ${\cal G}^{f}$ being the generator of a given semigroup ${\cal T}^{f}$.
To simplify the problem {\em we will assume that each ${\cal T}^{f}$ is the semigroup asociated to a transition system ${\cal Q}^{f}$}. The consistency conditions that have to be imposed on the collection $\{{\cal Q}^{f}; f \in {\cal S}\}$ so that it leads to a set-indexed transition system ${\cal Q}$ (and hence to a ${\cal Q}$-Markov process) are given in Section 4 (Assumptions 2 and 3). Using an integral form of the Kolomogorov-Feller equations, we prove that these conditions can be expressed equivalently in terms of the generators ${\cal G}^{f}$ and the semigroups ${\cal T}^{f}$. ({\em We conjecture that the same formulas will be valid in the general case when there is no transition system ${\cal Q}^{f}$ associated to the semigroup ${\cal T}^{f}$.}) 

Throughout, we shall illustrate our results with three examples of set-Markov processes: processes with independent increments, the empirical process, and the Dirichlet process.

We note in passing that set-Markov processes satisfy a type of strong Markov property which is studied in a separate paper \cite{balan01}.

\section{The Set-Indexed Framework}

This section introduces the general definitions, properties and
assumptions that are used in the theory of set-indexed
processes, as presented in \cite{ivanoff-merzbach00}. We will also 
give several examples, in addition to the multiparameter case.

Let ${\cal A}$ be a {\bf semilattice} of closed subsets of a compact
Hausdorff topological space $T$ (i.e. ${\cal A}$ is closed under
arbitrary  intersections), which contains the empty set and the
space $T$ itself, but does not contain disjoint (non-empty)
sets. In addition, we assume that the collection ${\cal A}$ is
{\bf separable from above}, in the sense that any set $A \in
{\cal A}$ can be approximated  from above as $$A=
\cap_{n}g_{n}(A); \ \ g_{n+1}(A) \subseteq g_{n}(A), A \subseteq
g_{n}(A)^{0} \ \forall n$$   where the approximation set
$g_{n}(A)$ can be written as a finite union of sets that lie in
a finite sub-semilattice ${\cal A}_{n}$ of ${\cal A}$; moreover,
${\cal A}_{n} \subseteq {\cal A}_{n+1} \ \forall n$ and $g_{n}$
preserves arbitrary intersections and finite unions i.e. 
$g_{n}(\cap_{\alpha \in \Lambda}A_{\alpha})= \cap_{\alpha \in
\Lambda}g_{n}(A_{\alpha}), \ \forall A_{\alpha} \linebreak \in
{\cal A}$; and  $\cup_{i=1}^{k}A_{i}= \cup_{j=1}^{m}A'_{j}
\Rightarrow
\cup_{i=1}^{k}g_{n}(A_{i})=\cup_{j=1}^{m}g_{n}(A'_{j}), \
\forall A_{i},A'_{j} \in {\cal A}$. By convention,
$g_{n}(\emptyset)= \emptyset$.

There are many examples of classes of sets which have these
properties.

\vspace{2mm}
{\bf Examples 1.}
\begin{enumerate}
\item $T=[0,1]^{d}, {\cal A}=\{[0,z];z \in T\} \cup \{\emptyset\}$.

\item $T=[0,1]^{d}, {\cal A}= \{A; A \ {\rm a \ compact \ lower \ 
layer \ in} \ T \} \cup \{\emptyset\}$. ($A$ is a lower layer if
 $z \in A \Rightarrow [0,z] \subseteq A$)

\item $T=[a,b]^{d}, a<0<b, {\cal A}=\{[0,z];z \in T\} \cup
\{\emptyset\}$.

\item $T=[a,b]^{d}, a<0<b, {\cal A}= \{A; A \ {\rm a \ compact \ lower
\  layer \ in} \ T \} \cup \{\emptyset\}$.  

\item $T= \overline{B(0,t_{0})}$ (compact ball in ${\bf R}^{3}$),
${\cal A}= \{A_{R,t}; R:=`[a,b] \times 
[c,d]', 0 \leq a < b < 2 \pi, - \pi \leq c < d \leq \pi,t \in [0,t_{0}]
\}$, where the set  
$$A_{R,t}:= \{(r \cos \theta \cos \tau,r \sin
\theta \cos \tau,r \sin \tau); \theta \in [a,b], \tau \in
[c,d], r \in [0,t] \}$$ can be interpreted as the history of the
region 
$$R:=`[a,b] \times [c,d]'= \{( \cos \theta \cos
\tau, \sin \theta \cos \tau, \sin \tau); \theta \in [a,b], \tau
\in [c,d] \}$$ 
of the Earth from the beginning until time $t$.
(Here $\theta$ represents the longitude of the generic point in the
region $R$, while $\tau$ is the latitude.) Hence, ${\cal
A}$ can be identified with the history of the world until time $t_{0}$. 
\end{enumerate} 
\vspace{2mm}

Let $\emptyset':= \cap_{A \in {\cal A} \verb2\2 \{\emptyset\}}A$ be the 
minimal set in ${\cal A}$ ($\emptyset' \not = \emptyset$). The role
played by $\emptyset'$ will be  similar to the role played by $0$ in the
classical theory.

We will consider the following classes of sets generated by ${\cal A}$:
\begin{itemize}
\item ${\cal A}(u)$ is the class of all finite unions of sets in
${\cal A}$
\item ${\cal C}$ is the class of all sets of the form $C=A \verb2\2 B$
with $A \in {\cal A}, B \in {\cal A}(u)$ 
\item ${\cal C}(u)$ is the class of all finite unions of sets in
${\cal C}$
\end{itemize}
 
Note that ${\cal A}(u)$ is closed under finite intersections or finite
unions, ${\cal C}$ is a semialgebra and ${\cal C}(u)$ is the algebra
generated by ${\cal C}$. The value $X_{C}=X_{A \cup B}-X_{B}$ of an
(additive) set-indexed process $(X_{A})_{A \in {\cal A}}$ over the set
$C=A \verb2\2 B= (A \cup B) \verb2\2 B$ will play the role of the
increment $X_{t}-X_{s},s <t$ of a one-dimensional process $(X_{t})_{t
\in [0,a]}$. 

Any set $B \in {\cal A}(u)$ admits at least one {\bf extremal
representation} of the form $B= \cup_{i=1}^{n}A_{i}, A_{i} \in
{\cal A}$ with $A_{i} \not \subseteq \cup_{j \not =i}A_{j} \ \forall i$.
For a set $C \in {\cal C}$ we will say that the representation
$C=A \verb2\2 B$ is extremal if the representation of $B$ is
extremal.

Since the finite sub-semilattices of the indexing
collection ${\cal A}$ play a very important role in the
theory of set-indexed processes, having an appropriate ordering on
the sets of a finite sub-semilattice proves to be a very useful
tool in handling these objects. The ordering that we have in
mind will be defined in such a manner that a set is never
numbered before any of its subsets. We will call such an
ordering {\bf consistent} (with the strong past). 
More precisely, if ${\cal A'}$ is a finite sub-semilattice of
${\cal A}$ we set $A_{0}:=\emptyset'$ and $A_{1}:= \cap_{A \in {\cal
A'}}A$ (note that the sets $A_{0}$ and $A_{1}$ are not necessarily
distinct). Proceeding inductively, assuming that the distinct sets
$A_{1}, \ldots, A_{i} \in {\cal A'}$ have already been counted, choose
$A_{i+1} \in {\cal A'} \verb2\2 \{A_{1} \ldots, A_{i}\}$ such
that if there exists a set $A \in {\cal A'}$ with $A \subseteq A_{i+1},
A \neq A_{i+1}$ then $A=A_{j}$ for some $j \leq i$. It is clear that
such an ordering always exists although in general, it is not unique. 

If $\{A_{0}=\emptyset',A_{1}, \ldots,A_{n}\}$ is a consistent ordering
of a finite sub-semilattice ${\cal A'}$, the set
$C_{i}=A_{i} \verb2\2 \cup_{j=0}^{i-1}A_{j} \in {\cal C}$ is called
{\bf the left neighbourhood} of the set $A_{i}$ (we make the convention
that the left neighbourhood of $A_{0}=\emptyset'$ is itself). The
definition of the left neighbourhood does not depend on the ordering
since one can show that $C_{i}=A_{i} \verb2\2 (\cup_{A \in {\cal A'},A_{i}
\not \subseteq A}A)$.

\vspace{2mm}
{\bf Comments 1.}
\begin{enumerate}
\item If $B_{1}, \ldots,B_{m} \in {\cal
A}(u)$ are such that $B_{1} \subseteq \cdots \subseteq B_{m}$
then there exists a finite sub-semilattice ${\cal
A'}$ of ${\cal A}$ and a consistent ordering $\{A_{0}=\emptyset',A_{1},
\ldots,A_{n}\}$ of ${\cal A'}$ such that
$B_{l}=\cup_{j=0}^{i_{l}}A_{j}$, for some $0< i_{1} \leq \ldots \leq
i_{m}=n$. 

\item If ${\cal A'},{\cal A''}$ are two
finite sub-semilattices of ${\cal A}$ such that ${\cal A'}
\subseteq {\cal A''}$, there exists a consistent ordering
$\{A_{0}=\emptyset', A_{1}, \ldots,A_{n}\}$ of ${\cal A''}$ such that
if ${\cal A'}=\{A_{i_{0}}=\emptyset',A_{i_{1}}, \ldots,A_{i_{m}}\}$
with $0=i_{0} <i_{1} \leq \ldots \leq i_{m}$, then
$\cup_{s=1}^{l}A_{i_{s}}=\cup_{j=1}^{i_{l}}A_{j}$ for any $l=1,
\ldots,m$. 
\end{enumerate}
\vspace{2mm}

Let us consider now a complete probability space $(\Omega, {\cal F}, P)$
and let $X:=(X_{A})_{A \in {\cal A}}$ be an ${\bf R}$-valued process
defined on this space, indexed by the class ${\cal A}$.
 
We say that the process $X$ {\bf has an (almost surely) unique
additive extension to ${\cal A}(u)$ } if whenever the set $B \in
{\cal A}(u)$ can be written as $B=
\cup_{i=1}^{n}A_{i}=\cup_{j=1}^{m}A'_{j}$ with
$A_{1}, \ldots,A_{n},A'_{1}, \ldots, A'_{m} \in {\cal A}$ we
have
\begin{eqnarray*}
\lefteqn{ \sum_{i=1}^{n}X_{A_{i}}- \sum_{1 \leq i_{1}
<i_{2} \leq n}X_{A_{i_{1}} \cap A_{i_{2}}}+ \cdots
+(-1)^{n+1}X_{A_{1} \cap \ldots \cap A_{n}} = } \\
 & & \sum_{j=1}^{m}X_{A'_{j}}- \sum_{1 \leq j_{1}
<j_{2} \leq m}X_{A'_{j_{1}} \cap A'_{j_{2}}}+ \cdots
+(-1)^{m+1}X_{A'_{1} \cap \ldots \cap A'_{m}} \  {\rm a.s.} 
\end{eqnarray*}
In this case, outside a set of measure $0$, we can define
$X_{B}$ as being either one of the members of the above
almost sure equality.  
In order to verify that a process has a unique additive
extension to ${\cal A}(u)$ it is enough to prove the previous almost
sure equality only in the case when $n=1$. The general case will follow
since if $\cup_{i=1}^{n}A_{i}= \cup_{j=1}^{m}A'_{j}$, then each of the
sets $A_{i_{1}} \cap A_{i_{2}} \cap \ldots \cap A_{i_{k}} \in {\cal A}$
can be written as the union $\cup_{j=1}^{m} (A_{i_{1}} \cap A_{i_{2}} 
\cap \ldots \cap A_{i_{k}} \cap A'_{j})$.

Similarly, the process $X$ is said to {\bf
have an (almost surely) unique additive extension to ${\cal C}$
} if whenever the set $C \in {\cal C}$ can be written as $C=A
\verb2\2 B= A' \verb2\2 B'$ with 
$A,A' \in {\cal A}$ and $B,B' \in {\cal A}(u)$ we have
$$X_{A}-X_{A \cap B}=X_{A'}-X_{A' \cap B'} \  {\rm a.s.}$$
The additive extension to ${\cal C}(u)$ is defined in the obvious
manner.

A set-indexed process $X:=(X_{A})_{A \in {\cal A}}$ with a unique
additive extension to ${\cal C}(u)$, for which $X_{C_{1}}, \ldots,
X_{C_{n}}$ are independent whenever the sets $C_{1}, \ldots,C_{n} \in
{\cal C}$ are disjoint, is called a {\bf process with independent
increments}.

An example of a process which has a unique additive extension to
${\cal C}(u)$ is the {\bf empirical process} of size $n$, corresponding
to a probability measure $F$ on $T$, defined by $X_{A}:=
\sum_{j=1}^{n}I_{ \{Z_{j} \in A \} },A \in {\cal A}$, where $(Z_{j})_{j
\geq 1}$ are i.i.d. $T$-valued random variables with common distribution
$F$.

Another example is the {\bf Dirichlet process} with parameter measure $\alpha$, where
$\alpha$ is a finite positive measure on $\sigma({\cal A})$ (see Ferguson \cite{ferguson73}).
This process is in fact almost surely {\em countably} additive on $\sigma({\cal A})$ and takes values in $[0,1]$; moreover $X_{T}=1$ a.s. Its
finite dimensional distribution of this process over any disjoint sets
$A_{1}, \ldots, A_{k}$ with $A_{i} \in \sigma({\cal A})$ is the (non-singular) Dirichlet distribution
with parameters $(\alpha(A_{1}), \ldots, \alpha(A_{k}); \alpha((\cup_{i=1}^{k}A_{i})^{c}))$
(see Ferguson \cite{ferguson73}).

In what follows we will examine the information structure
that can be associated with a set-indexed process.

A collection $({\cal F}_{A})_{A \in
{\cal A}}$ of sub-$\sigma$-fields of ${\cal F}$ is called a {\bf
filtration} if ${\cal F}_{A} \subseteq {\cal F}_{A'}$ whenever
$A,A' \in {\cal A}, A \subseteq A'$. We will consider
only complete filtrations. 
An ${\cal A}$-indexed filtration $({\cal F}_{A})_{A \in
{\cal A}}$   can be extended   to a filtration indexed by ${\cal
A}(u)$ by defining for each $B \in {\cal A}(u)$,
\begin{equation}
{\cal F}_{B}:=\bigvee_{A \in {\cal A}, A \subseteq B}{\cal F}_{A} 
\end{equation}
A set-indexed process $X$ is {\bf adapted} with respect to the
filtration $({\cal F}_{A})_{A \in {\cal A}}$ if $X_{A}$ is ${\cal
F}_{A}$ -measurable  for any $A \in {\cal A}$. 
If $X$ is adapted and has a unique additive extension to ${\cal A}(u)$
then $X_{B}$ is ${\cal F}_{B}$-measurable for any $B \in {\cal A}(u)$.
Given a set-indexed process $X$, the {\bf minimal filtration} with
respect to which $X$ is adapted is given by
${\cal F}_{B}:=\sigma(\{X_{A}; A \in {\cal A}, A \subseteq B \})$.

Any map $f: [0,a] \rightarrow {\cal A}(u)$ which is increasing 
with respect to the partial order induced by the set-inclusion is
called a {\bf flow}. 

\vspace{2mm}
{\bf Definition 1.}
{\em A flow $f: [0,a] \rightarrow {\cal A}(u)$ is
\begin{description}
\item[a)] {\bf continuous} if for any $t \in [0,a]$ and for any
decreasing sequence $(t_{n})_{n}$ with $\lim_{n \rightarrow
\infty}t_{n}=t$ we have $f(t)= \cap_{n}f(t_{n})$, and for any
increasing sequence $(t_{n})_{n}$ with $\lim_{n \rightarrow
\infty}t_{n}=t$ we have $f(t)= \overline{\cup_{n}f(t_{n})}$. 
\item[b)] {\bf simple} if it is continuous and there exists a
partition $0=t_{0}<t_{1}< \ldots <t_{n}=a$ and flows
$f_{i+1}:[t_{i}, t_{i+1}] \rightarrow {\cal A}$, $i=0, \ldots n-1$ such
that $f(0)=\emptyset'$ and $f(t)=\cup_{j=1}^{i}f_{j}(t_{j}) \cup
f_{i+1}(t), t \in [t_{i},t_{i+1}],i=0, \ldots,n-1$. (In other words, a
simple flow is piecewise ${\cal A}$-valued.) 
\end{description} }
\vspace{2mm}

If ${\cal A'}$ is a finite sub-semilattice and ord$=\{\emptyset'=A_{0},
A_{1}, \ldots, A_{n}\}$ is a consistent ordering of ${\cal A'}$, we say
that a simple flow $f$ {\bf connects the sets of the semilattice}
${\cal A'}$, in the sense of the ordering ord, if $f(t)=
\cup_{j=1}^{i}A_{j} \cup f_{i+1}(t), t \in [t_{i},t_{i+1}]$, where
$t_{0}=0< t_{1} < \ldots < t_{n}=a$ is a partition of the domain of
definition of $f$ and $f_{i+1}:[t_{i},t_{i+1}] \rightarrow {\cal A}$ are
continuous flows with $f_{i+1}(t_{i})=A_{i}, f_{i+1}(t_{i+1})=A_{i+1}$.

\vspace{2mm}
{\bf Lemma 1.} {\em (Lemma 5.1.7, \cite{ivanoff-merzbach00})
For every finite sub-semilattice ${\cal A}'$ and for each
consistent ordering ord of ${\cal A'}$, there exists a simple flow $f$
which connects the sets of the semilattice ${\cal A'}$, in the sense of
the ordering ord.}
\vspace{2mm}

{\bf Comment 2.}
A consequence of the previous lemma is the following: given $B_{1},
\ldots, B_{m} \in {\cal A}(u)$ such that $B_{1} \subseteq \ldots
\subseteq B_{m}$ there exists a simple flow $f$ and   $t_{1}\leq \ldots
\leq t_{m}$ such that $B_{i}=f(t_{i});i=1, \ldots,m$.

\section{Set-Markov Processes}

In this section we will introduce the definition of the set-Markov 
property. Two immediate consequences of the definition will be: 
(a) any process with independent increments is set-Markov; and 
(b) the set-Markov property is equivalent to the classical Markov
property on every flow. Finally we will prove that any set-Markov
process is also sharp Markov.

If ${\cal F}, {\cal G}, {\cal H}$ are three sub-$\sigma$-fields of
the same probabilistic space, we will use the notation ${\cal F} \perp
{\cal H} \ | \ {\cal G}$ if ${\cal F}$ and ${\cal H}$ are
conditionally independent given ${\cal G}$. 

\vspace{2mm}
{\bf Definition 2.}
{\em Let $X:=(X_{A})_{A \in {\cal A}}$ be a set-indexed process with a
unique additive extension to ${\cal C}(u)$ and $({\cal F}_{A})_{A \in
{\cal A}}$ its minimal filtration.
We say that the process $X$ is {\bf set-Markov} if 
$\forall A \in {\cal A}, \ \forall B \in {\cal A}(u)$,
${\cal F}_{B} \perp \sigma(X_{A \verb2\2 B}) \ | \ \sigma(X_{B})$.}
\vspace{2mm}

\noindent It is easy to see that in the classical case the set-Markov
property is equivalent to the usual Markov property.

\vspace{2mm}
{\bf Examples 2.}
\begin{enumerate}
\item Any process $X:=(X_{A})_{A \in {\cal A}}$ with independent
increments is set-Markov since $\forall A \in {\cal A}, \ \forall B \in {\cal A}(u)$, $X_{A \verb2\2 B}$ is independent of ${\cal F}_{B}$.

\item The empirical process of size $n$ (corresponding to $F$) is set-Markov since $\forall A \in {\cal A}, \ \forall B \in {\cal A}(u)$, for any partition $B= \cup_{i=1}^{p}C_{i}, C_{i} \in {\cal C}$ and for any $l,k_{1}, \ldots,k_{p} \in \{0,1, \ldots, n\}; k:=\sum_{i=1}^{p}k_{i}$ with $k+l \leq n$
$$P[X_{A \verb2\2 B}={{l} \over {n}} | X_{C_{i}}={{k_{i}} \over
{n}};i=1, \ldots,p]= P[X_{A \verb2\2 B}={{l} \over {n}} | X_{B}={{k} \over {n}}]$$   
both sides being equal to the value at $l$ of the binomial distribution with $n-k$ trials and ${{F(A \verb2\2 B)} \over {1-F(B)}}$ probability of success.

\item The Dirichlet process $X:=(X_{A})_{A \in \sigma({\cal A})}$ with parameter measure $\alpha$ is set-Markov since
$\forall A \in {\cal A}, \ \forall B \in {\cal A}(u)$, for any partition $B= \cup_{i=1}^{p}C_{i}, C_{i} \in {\cal C}$ and for any $y,x_{1}, \ldots,x_{p} \in [0,1]; x:=\sum_{i=1}^{p}x_{i}$ with $x+y \leq 1$
$$P[X_{A \verb2\2 B} \leq y|X_{C_{i}}=x_{i}; i=1, \ldots,p]=P[X_{A \verb2\2 B} \leq y|X_{B}=x]$$
both sides being equal to the value at ${{y} \over {1-x}}$
of the Beta distribution with parameters $(\alpha(A \verb2\2 B); \alpha((A \cup B)^{c}))$
(we use property {\bf 7.7.3}, p. 180, \cite{wilks63} of the Dirichlet distribution).

\end{enumerate}

We shall make repeated use of the following elementary result.

\vspace{2mm}
{\bf Lemma 2.}
{\em Let ${\cal G'} \subseteq {\cal G}$ be two $\sigma$-fields in the
same probability space and $X, Y$ two random vectors on this space such
that $Y$ is ${\cal G}$-measurable. Suppose that $E[f(X)|{\cal
G}]=E[f(X)|{\cal G}']$ for every bounded measurable function $f$. Then
$E[h(X,Y)|{\cal G}]=E[h(X,Y)|{\cal G}',Y]$ for every bounded measurable
function $h$. } 
\vspace{2mm}

{\bf Proposition 1.}
{\em Let $X:=(X_{A})_{A \in {\cal A}}$ be a set-indexed
process with a unique additive extension to ${\cal C}(u)$ and
$({\cal F}_{A})_{A \in {\cal A}}$ its minimal filtration. The process 
$X$ is set-Markov if and only if $\forall B, B' \in {\cal A}(u), B
\subseteq B'$,  
${\cal F}_{B} \perp \sigma(X_{B'}) \ | \ \sigma(X_{B})$.}
\vspace{2mm}

\noindent{\bf Proof}: Using Lemma 2 and the
additivity of the process, it follows that $X$ is set-Markov if and
only if
$\forall A
\in {\cal A} \ \forall B \in {\cal A}(u)$,  
${\cal F}_{B} \perp \sigma(X_{A \cup B}) \ | \ \sigma(X_{B})$. 
(Write $X_{A \cup B}=X_{A \verb2\2 B} + X_{B}$ and use the fact that
$X_{B}$ is ${\cal F}_{B}$-measurable.) 
Let $B,B' \in {\cal A}(u)$ be such that $B \subseteq B'$. Say
$B'= \cup_{i=1}^{k}A_{i}, A_{i} \in {\cal A}$ is an arbitrary
representation. Then $B'=B \cup B'= B \cup \cup_{i=1}^{k}
A_{i}$ and the result follows by induction on $k$.
$\Box$

\vspace{2mm}
{\bf Comment 3.} Using Proposition 1 and Lemma 2, we can say that a
process $X:=(X_{A})_{A \in {\cal A}}$ is set-Markov if and only if 
$\forall B, B' \in {\cal A}(u), B \subseteq B'$,  
${\cal F}_{B} \perp \sigma(X_{B' \verb2\2 B}) \ | \ \sigma(X_{B})$.
\vspace{2mm}

The following result says that a process is set-Markov if and only if
it is Markov (in the usual sense) along any flow. Moreover, it suffices
to restrict our attention only to simple flows.

\vspace{2mm}
{\bf Proposition 2.}
{\em Let $X:=(X_{A})_{A \in {\cal A}}$ be a set-indexed
process with a unique additive extension to ${\cal C}(u)$ and 
$({\cal F}_{A})_{A \in {\cal A}}$ its minimal filtration. 
Then $X$ is set-Markov if and only if for every simple flow
$f:[0,a] \rightarrow {\cal A}(u)$ the process
$X^{f}:=(X_{f(t)})_{t \in [0,a]}$ is Markov with respect to the
filtration $({\cal F}_{f(t)})_{t \in [0,a]}$.
(For  necessity, we can consider any flow, not only the
simple ones.)}
\vspace{2mm}

\noindent{\bf Proof}: The process $X^{f}$ is Markov with respect to the
filtration $({\cal F}_{f(t)})_{t}$ if and only if $\forall s,t \in
[0,a],s <t, {\cal F}_{f(s)} \perp \sigma(X_{f(t)}) \ | \
\sigma(X_{f(s)})$.  (We note  that the filtration $({\cal F}_{f(t)})_{t
\in [0,a]}$ is not the minimal filtration associated to the process
$X^{f}$.) This is equivalent to the set-Markov property since
we know that whenever the sets $B,B' \in {\cal A}(u)$ are such that $B
\subseteq B'$ there exists a simple flow $f$ and some $s <t$ such that
$f(s)=B$ and $f(t)=B'$ (Comment 2). $\Box$

The preceding proposition, while simple, is crucial since (as noted in
the Introduction)  it provides us with the means to define the generator
of a set-Markov process. This will be done in Section \ref{Generator}.

We now show that every set-Markov process satisfies the sharp-Markov
property defined in ~\cite{ivanoff-merzbach98}. In analogy with the
minimal filtration $({\cal F}_{B})_{B \in {\cal A}(u)}$ we define the
following $\sigma$-fields, for an arbitrary set
$B \in {\cal A}(u)$: 
$${\cal F}_{\partial B}:=\sigma(\{X_{A}; A \in {\cal A}, 
A \subseteq B, A \not \subseteq B^{0} \})$$
$${\cal F}_{B^{c}}:=\sigma(\{X_{A}; A \in {\cal A}, A \not \subseteq B
\})$$

\vspace{2mm}
{\bf Definition 3.}
{\em Let $X:=(X_{A})_{A \in {\cal A}}$ be a set-indexed process with a
unique additive extension to ${\cal C}(u)$ and $({\cal F}_{A})_{A \in
{\cal A}}$ its minimal filtration.
We say that the process $X$ is {\bf sharp-Markov} if 
$\forall B \in {\cal A}(u)$,   
${\cal F}_{B} \perp {\cal F}_{B^{c}} \ | \ {\cal F}_{\partial B}$. }
\vspace{2mm}

The next lemma will be essential in proving that a set-Markov process is sharp Markov.

\vspace{2mm}
{\bf Lemma 3.}
{\em If $X:=(X_{A})_{A \in {\cal A}}$ is a set-Markov process, then
$\forall B \in {\cal A}(u) \ \forall A_{i} \linebreak \in {\cal A};i=1,
\ldots,n$, ${\cal F}_{B} \perp \sigma(X_{A_{1} \verb2\2 B}, \ldots,
X_{A_{n} \verb2\2 B}) \ | \ \sigma(X_{B})$. }
\vspace{2mm}

\noindent{\bf Proof}: Let $B \in {\cal A}(u)$,$A_{1}, \ldots A_{n} \in
{\cal A}$ and $h:{\bf R}^{n} \rightarrow {\bf R}$ an arbitrary bounded
measurable function. Without loss of generality we can assume
that $A_{i} \not \subseteq B \  \forall i$; say
$B=\cup_{j=n+1}^{m}A_{j}$, $A_{j} \in {\cal A}$. Let ${\cal
A'}$ be the smallest finite sub-semilattice which
contains $A_{1}, \ldots A_{m}$, $\{A'_{0}=\emptyset',A'_{1},
\ldots, A'_{p}\}$ a consistent ordering of ${\cal A'}$, and
$C'_{i}$ the left neighborhood of $A'_{i}$ in ${\cal A'}$. 
Say $A_{j}=A'_{i_{j}}$ for $j=1 \ldots n$; then $A_{j} \verb2\2
B=A'_{i_{j}} \verb2\2 B=\cup_{i \in I_{j}} C'_{i}$ with $I_{j}
\subseteq \{1, \ldots, i_{j}\}$ and $X_{A_{j} \verb2\2
B}=\sum_{i \in I_{j}}X_{C'_{i}}$.
 
Therefore we can say that $h(X_{A_{1} \verb2\2 B}, \ldots,
 X_{A_{n} \verb2\2 B})=h_{1}(X_{C'_{l_{1}}}, \ldots,
X_{C'_{l_{s}}})$, for a certain bounded measurable function
$h_{1}$ and some $l_{1} \leq \ldots \leq l_{s},C'_{l_{i}} \not
\subseteq B$. In order to simplify the notation, let us denote
$D_{i}:=C'_{l_{i}}, i=1, \ldots,s$.  
 Let $B_{i}=B \cup \cup_{k=1}^{i}D_{k}$ for $i=1, \ldots , s$.
 Then $D_{i}=B_{i} \verb2\2 B_{i-1}$ and
 $X_{D_{i}}=X_{B_{i}}-X_{B_{i-1}}$; hence  $X_{D_{1}}, \ldots,
X_{D_{s-1}}$ are ${\cal F}_{B_{l_{s-1}}}$-measurable.
Because each $D_{i}$ is the left neighbourhood of $A'_{l_{i}}$,
we also have $B_{i}=B \cup
\cup_{k=1}^{i}A'_{l_{k}}$ and therefore $D_{i}=(A'_{l_{i}} \cup
B_{i-1}) \verb2\2 B_{i-1}=A'_{l_{i}} \verb2\2 B_{i-1}$.  
Using the equivalent definition of the set-Markov property given by
Comment 3, we have  $E[f(X_{D_{s}})|{\cal
F}_{B_{s-1}}]=E[f(X_{D_{s}})|X_{B_{s-1}}]$ for any bounded
measurable function $f$. Now we are in the position to apply
Lemma 2 to get   
$$E[h_{1}(X_{D_{1}}, \ldots,X_{D_{s}})|{\cal
F}_{B_{s-1}}]=E[h_{1}(X_{D_{1}}, \ldots, X_{D_{s}})|X_{B_{s-1}},
X_{D_{1}}, \ldots, X_{D_{s-1}}].$$
Hence
$E[h(X_{A_{1} \verb2\2 B}, \ldots, X_{A_{n} \verb2\2
B})|{\cal F}_{B}] = E[E[h_{1}(X_{D_{1}}, \ldots, X_{D_{s}})|{\cal
F}_{B_{s-1}}]|{\cal F}_{B}] = \linebreak E[E[h_{1}(X_{D_{1}}, \ldots,
X_{D_{s}})|X_{B_{s-1}},X_{D_{1}}, \ldots, X_{D_{s-1}}]|{\cal F}_{B}]$.
Writing $X_{B_{s-1}}=X_{B}+ \sum_{k=1}^{s-1}X_{D_{k}}$
we get $E[h(X_{A_{1} \verb2\2 B}, \ldots, X_{A_{k} \verb2\2
B})|{\cal F}_{B}]=E[h_{2}(X_{D_{1}}, \ldots,
X_{D_{s-1}},X_{B})| \linebreak {\cal F}_{B}]$.
Continuing in the same manner, reducing at each step another
set $D_{i}$, we finally get $E[h(X_{A_{1} \verb2\2 B}, \ldots,
X_{A_{n} \verb2\2 B})|{\cal F}_{B}]=E[h(X_{A_{1} \verb2\2 B},
\ldots, X_{A_{n} \verb2\2 B})|X_{B}]$. $\Box$

\vspace{2mm}
{\bf Proposition 3.}
{\em Any set-Markov process is sharp Markov.}
\vspace{2mm}

\noindent{\bf Proof:} Let $B \in {\cal A}(u), A_{i} \in {\cal A},A_{i}
\not \subseteq B; i=1, \ldots,n$ and $h:{\bf R}^{n} \rightarrow {\bf R}$
an arbitrary bounded measurable function. Writing $X_{A_{i}}=X_{A_{i}
\cap B} +X_{A_{i} \verb2\2 B}$ and using Lemma 2 we
can say that $E[h(X_{A_{1}}, \ldots, X_{A_{n}})|{\cal
F}_{B}]=E[h(X_{A_{1}}, \ldots,  X_{A_{n}})|X_{B}, \newline
X_{A_{i} \cap B}, i=1, \ldots n]$, which is ${\cal F}_{\partial
B}$-measurable, by Lemma 2.4, \cite{ivanoff-merzbach98}.
$\Box$ 

In what follows we will give an important characterization of the
set-Markov processes that will be instrumental for the
construction of these processes.  

We will need the following elementary result.

\vspace{2mm}
{\bf Lemma 4.}
{\em Let ${\cal F}_{i},i=1, \ldots,4$ be four $\sigma$-fields
in the same probability space. If
${\cal F}_{1} \perp {\cal F}_{3} \  | \  {\cal F}_{2}$ and
${\cal F}_{1} \vee {\cal F}_{2} \perp {\cal F}_{4} \  | \ 
{\cal F}_{3}$, then ${\cal F}_{1} \perp {\cal F}_{3} \vee
{\cal F}_{4} \  | \  {\cal F}_{2}$. }
\vspace{2mm}

{\bf Proposition 4.}
{\em A set-indexed process $X:=(X_{A})_{A \in {\cal A}}$ with a unique
additive extension to ${\cal C}(u)$, is set-Markov if and only
if for every finite sub-semilattice ${\cal A'}$ and for every
consistent ordering $\{A_{0}=\emptyset',A_{1}, \ldots, A_{n}\}$ of
${\cal A'}$ }
\begin{equation}
\label{one-term}
\sigma(X_{A_{0}},X_{A_{0} \cup A_{1}}, \ldots,
X_{\cup_{j=0}^{i-1}A_{j}}) \perp
\sigma(X_{\cup_{j=0}^{i+1}A_{j}})
|\sigma(X_{\cup_{j=0}^{i}A_{j}}) \ \forall i=1, \ldots, n-1
\end{equation} 

\noindent{\bf Proof}: We will use the characterization of the set-Markov
property given by Proposition 1.
Necessity follows immediately.

For sufficiency note first that equation (\ref{one-term}) implies a
similar equation where the union of the first $i+1$  sets is
replaced by the union of the first $i+p$  sets. In fact it is
easily shown by induction on $p$ and using Lemma 4, that   
\begin{equation}
\label{CO-pterms}
\sigma(X_{A_{0}},X_{A_{0} \cup A_{1}}, \ldots,
X_{\cup_{j=0}^{i-1}A_{j}}) \perp
\sigma(X_{\cup_{j=0}^{i+k}A_{j}};k=1, \ldots,p)
 | \sigma(X_{\cup_{j=0}^{i}A_{j}}) \  \forall i
\end{equation}  

Consider now arbitrary sets $B,B' \in {\cal A}(u)$ with $B
\subseteq B'$. By a monotone class argument it is enough to show that
$\forall A'_{l} \in {\cal A}, A'_{l} \subseteq B; l=1, \ldots,m$,
$\sigma(X_{A'_{1}}, \ldots, X_{A'_{m}}) \perp \sigma(X_{B'})
|\sigma(X_{B})$.  Without loss of generality we may assume that
$\{A'_{0}=\emptyset',A'_{1}, \ldots, A'_{m}\}$ is a finite
sub-semilattice and that the ordering
$\{A'_{0}=\emptyset',A'_{1}, \ldots, A'_{m}\}$ is consistent. Because
the process $X$ has a unique additive extension to ${\cal A}(u)$, there
exists a bijective map $\psi$ such that $(X_{A'_{1}},X_{A'_{2}}, \ldots,
X_{A'_{m}})=\psi(X_{A'_{1}},X_{A'_{1} \cup A'_{2}}, \ldots,
X_{\cup_{l=1}^{m}A'_{l}})$ a.s.. Consequently, we have
$\sigma(X_{A'_{1}},X_{A'_{2}}, \ldots, X_{A'_{m}})=\sigma(X_{A'_{1}},
X_{A'_{1} \cup A'_{2}}, \ldots, X_{\cup_{l=1}^{m}A'_{l}})$.
  
By Comment 1.1, there exists a finite sub-semilattice ${\cal
A}'$ of ${\cal A}$ and a consistent ordering $\{A_{0}= \emptyset',A_{1},
\ldots, A_{n}\}$ of ${\cal A'}$, such
that $A'_{1}=\cup_{j=0}^{i_{1}}A_{j},A'_{1} \cup
A'_{2}=\cup_{j=0}^{i_{2}}A_{j}, \ldots,
\cup_{l=1}^{m}A'_{l}=\cup_{j=0}^{i_{m}}A_{j}, 
B=\cup_{j=0}^{i_{m+1}}A_{j},B'=\cup_{j=0}^{i_{m+2}}A_{j}$ for
some $i_{1} \leq i_{2} \leq \ldots \leq i_{m+2}$. Using 
(\ref{CO-pterms}), it follows that
$\sigma(X_{A'_{1}}, X_{A'_{1} \cup A'_{2}}, \ldots,
X_{\cup_{l=1}^{m}A'_{l}}) \perp \sigma(X_{B'}) | \sigma(X_{B})$. $\Box$

\section{Construction of the Process}

In this section we will introduce a special class of set-Markov
processes for which the mechanism of transition from one
state to another is completely known, and we will construct such a
process. In light of Proposition 2, we will also
determine the necessary and sufficient conditions that have to be
imposed on a family of one-dimensional transition systems, indexed by
a collection of simple flows, such that on each simple flow $f$ from
the chosen collection, the corresponding Markov process has the law
of $X^f$, where $X$ is set-Markov (i.e., under what circumstances a
class of one dimensional transition systems determines a set-Markov
process).

Let ${\cal B}({\bf R})$ denote the Borel subsets of ${\bf R}$.
We begin with the definition of the transition system.

{\bf Definition 4.}
{\em
\begin{description}
\item[(a)]
For each $B,B' \in {\cal A}(u),B \subseteq B'$ let
$Q_{BB'}(x; \Gamma), x \in {\bf R}, \Gamma \in {\cal B}({\bf R})$, be
a transition probability on ${\bf R}$ i.e., $Q_{BB'}(x;
\cdot)$ is a probability measure $\forall x$, and $Q_{BB'}( \cdot ;
\Gamma)$ is measurable $\forall \Gamma$. The family ${\cal
Q}:=(Q_{BB'})$ of all these transition probabilities is called a
{\bf transition system} if $\forall B \in {\cal A}(u), Q_{BB}(x;
\cdot)= \delta_{x}$ and $\forall B,B',B'' \in {\cal A}(u),B \subseteq
B' \subseteq B''$   $$Q_{BB''}(x; \Gamma) = \int_{\bf
R}Q_{B'B''}(y;\Gamma)Q_{BB'}(x;dy) \ \forall x \in {\bf R}, 
\forall \Gamma \in {\cal B}({\bf R})$$

\item[(b)] Let ${\cal Q}:=(Q_{BB'})_{B,B' \in {\cal A}(u);  B
\subseteq B'}$ be a transition system. A set-indexed
process $X:=(X_{A})_{A \in {\cal A}}$ with a
unique  additive extension to ${\cal C}(u)$ is called  {\bf ${\cal
Q}$-Markov} if $ \forall B,B' \in {\cal A}(u), B \subseteq B'$
$$P[X_{B'} \in \Gamma|{\cal F}_{B}]=Q_{BB'}(X_{B};\Gamma) \ \forall
\Gamma \in {\cal B}({\bf R})$$   
where $({\cal F}_{A})_{A \in {\cal A}}$ is the minimal filtration of
the process $X$. 
\end{description}   }  
\vspace{2mm}

In other words, a ${\cal Q}$-Markov process is a set-Markov process
for which $Q_{BB'}$ is a version of the conditional distribution of
$X_{B'}$ given $X_{B}$, for every $B,B' \in {\cal A}(u),B
\subseteq B'$.

\vspace{2mm}
{\bf Examples 3.}
\begin{enumerate}
\item Any process $X:=(X_{A})_{A \in {\cal A}}$ with independent
increments is ${\cal Q}$-Markov with $Q_{BB'}(x; \Gamma):=F_{B'
\verb2\2 B}(\Gamma -x)$, where $F_{C}$ is the distribution of $X_{C},
C \in {\cal C}(u)$. The Poisson process and the Brownian motion are
the particular cases for which $F_{C}$ is a Poisson distribution
with mean $\Lambda_{C}$, respectively, a  normal distribution  with 
mean 0 and variance $\Lambda_{C}$.  

\item The empirical process of size $n$, corresponding to a
probability measure $F$, is ${\cal Q}$-Markov with 
$Q_{BB'}({{k} \over {n}}; \{ {{m} \over {n}} \}); k,m \in \{0,1, \ldots,n \}, k \leq m$ given by the value at $m-k$ of the binomial distribution with
$n-k$ trials and ${{F(B' \verb2\2 B)} \over {1-F(B)}}$ probability of
success. 

\item The Dirichlet process with parameter measure $\alpha$ is ${\cal Q}$-Markov with 
\linebreak $Q_{BB'}(x; [0,z]); x,z \in [0,1]$ given by the value at ${{y-x} \over {1-x}}$ of the Beta distribution with parameters $(\alpha(B' \verb2\2 B); \alpha(B'^{c}))$.
 
\end{enumerate} 
\vspace{2mm}

\noindent The following result gives equivalent characterizations of a ${\cal
Q}$-Markov process.

\vspace{2mm}
{\bf Proposition 5.}
{\em Let ${\cal Q}:=(Q_{BB'})_{B,B' \in {\cal A}(u); 
B \subseteq B'}$ be a transition system,  
$X:=(X_{A})_{A \in {\cal A}}$ a set-indexed process with a unique
additive extension to ${\cal C}(u)$ and initial distribution $\mu$, 
and $({\cal F}_{A})_{A \in {\cal A}}$ the minimal filtration of the
process $X$. The following statements are equivalent: 
\begin{description} 
\item[(a)] The process $X$ is ${\cal Q}$-Markov. 

\item[(b)] For every simple flow $f:[0,a] \rightarrow {\cal A}(u)$
the process $X^{f}:=(X_{f(t)})_{t \in [0,a]}$ is ${\cal
Q}^{f}$-Markov (with respect to $({\cal F}_{f(t)})_{t
\in [0,a]}$), where $Q_{st}^{f}:=Q_{f(s),f(t)}$. 

\item[(c)] For every finite sub-semilattice ${\cal A}'$, for every
consistent ordering $\{A_{0}=\emptyset', A_{1}, \ldots, A_{n}\}$ of
${\cal A'}$ and for every $i=1, \ldots,n-1$
$$\sigma(X_{A_{0}},X_{A_{1}}, X_{A_{1} \cup A_{2}}, \ldots,
X_{\cup_{j=1}^{i-1}A_{j}}) \perp \sigma(X_{\cup_{j=1}^{i+1}A_{j}}) \ |
\ \sigma(X_{\cup_{j=1}^{i} A_{j}})$$ and 
$Q_{\cup_{j=1}^{i}A{j},\cup_{j=1}^{i+1}A_{j}}$ is a version of the
conditional distribution of $X_{\cup_{j=1}^{i+1}A_{j}}$ given 
$X_{\cup_{j=1}^{i}A_{j}}$.

\item[(d)] For every finite sub-semilattice ${\cal A}'$ and for every
consistent ordering $\{A_{0}=\emptyset', A_{1}, \ldots, A_{n}\}$ of
${\cal A'}$  
$$P(X_{A_{0}} \in \Gamma_{0},X_{A_{1}} \in
\Gamma_{1}, X_{A_{1} \cup A_{2}} \in \Gamma_{2}, \ldots,
X_{\cup_{j=1}^{n}A_{j}} \in \Gamma_{n})=$$ 
$$\int_{{\bf R}^{n+1}}
I_{\Gamma_{0}} (x_{0}) \prod_{i=1}^{n}I_{\Gamma_{i}}(x_{i})
 Q_{\cup_{j=1}^{n-1}A_{j},
\cup_{j=1}^{n}A_{j}}(x_{n-1};dx_{n}) \ldots$$
$$Q_{A_{1},A_{1} \cup A_{2}}(x_{1};dx_{2})
Q_{\emptyset'A_{1}}(x_{0};dx_{1}) \mu(dx_{0})$$ 
for every $\Gamma_{0}, \Gamma_{1}, \ldots, \Gamma_{n} \in {\cal
B}({\bf R})$.

\item[(e)] For every finite sub-semilattice ${\cal A}'$, for every
consistent ordering $\{A_{0}=\emptyset', A_{1}, \ldots, A_{n}\}$ of
${\cal A'}$, if we denote with $C_{i}$ the left neighbourhood of
the set $A_{i}$, then
$$P(X_{C_{0}} \in \Gamma_{0},X_{C_{1}} \in
\Gamma_{1}, X_{C_{2}} \in \Gamma_{2}, \ldots, X_{C_{n}} \in
\Gamma_{n})=$$  $$\int_{{\bf R}^{n+1}} I_{\Gamma_{0}}(x_{0})
I_{\Gamma_{1}}(x_{1})
\prod_{i=2}^{n}I_{\Gamma_{i}}(x_{i}-x_{i-1})Q_{\cup_{j=1}^{n-1}A_{j},
\cup_{j=1}^{n}A_{j}}(x_{n-1};dx_{n}) \ldots$$ $$Q_{A_{1},A_{1} \cup
A_{2}}(x_{1};dx_{2}) Q_{\emptyset'A_{1}}(x_{0};dx_{1}) \mu(dx_{0})$$
for every $\Gamma_{0}, \Gamma_{1}, \ldots, \Gamma_{n} \in {\cal
B}({\bf R})$. 
\end{description} }
\vspace{2mm}

\noindent{\bf Proof}: The equivalences (a)-(b), (a)-(c) follow by 
arguments similar to those used to prove Proposition 2, respectively,
Proposition 4. The equivalence (c)-(d) follows
exactly as in the classical case. Finally, the equivalence (d)-(e)
follows by a change of variables, since $X_{C_{i}}=
X_{\cup_{j=1}^{i}A_{j}}- X_{\cup_{j=1}^{i-1}A_{j}}$ a.s. $\Box$

The general construction of a ${\cal Q}$-Markov process will be made
using increments: i.e., the sets in ${\cal C}$. The following 
assumption is necessary. It requires that the distribution of
the process over the left-neighbourhoods $C_{0},C_{1}, \ldots,C_{n}$
of a finite sub-semilattice, does not depend on the consistent
ordering of the semilattice.

\vspace{2mm}
{\bf Assumption 1.}
{\em If $\{A_{0}=\emptyset',A_{1}, \ldots,A_{n}\}$ and 
$\{A_{0}=\emptyset',A'_{1}, \ldots, A'_{n}\}$ are two consistent
orderings of the same finite sub-semilattice ${\cal A'}$ and $\pi$
is the permutation of $\{1, \ldots,n\}$ with $\pi(1)=1$ such that
$A_{i}=A'_{\pi(i)} \ \forall i$, then  
$$\int_{{\bf R}^{n+1}} I_{\Gamma_{0}}(x_{0})
I_{\Gamma_{1}}(x_{1})\prod_{i=2}^{n} I_{\Gamma_{i}}(x_{i}-x_{i-1})
Q_{\cup_{j=1}^{n-1}A_{j}, \cup_{j=1}^{n}A_{j}}(x_{n-1};dx_{n})
\ldots$$ $$Q_{A_{1},A_{1} \cup A_{2}}(x_{1};dx_{2})
Q_{\emptyset'A_{1}}(x_{0};dx_{1}) \mu(dx_{0})=$$  
$$\int_{{\bf R}^{n+1}} I_{\Gamma_{0}}(y_{0})I_{\Gamma_{1}}(y_{1})
\prod_{i=2}^{n}I_{\Gamma_{i}}(y_{\pi(i)}-y_{\pi(i)-1})Q_{\cup_{j=1}^{n-1}A'_{j},
\cup_{j=1}^{n}A'_{j}}(y_{n-1};dy_{n}) \ldots$$
$$Q_{A'_{1},A'_{1} \cup A'_{2}}(y_{1};dy_{2})
Q_{\emptyset'A'_{1}}(y_{0};dy_{1}) \mu(dy_{0})$$ 
for every $\Gamma_{0}, \Gamma_{1}, \ldots, \Gamma_{n} \in {\cal
B}({\bf R})$. } 
\vspace{2mm}

The finite dimensional distributions of a ${\cal Q}$-Markov process over the sets in ${\cal C}$ have to be defined so that they ensure the (almost sure) additivity of the process. The next result gives the definition of the finite dimensional distribution (of an additive ${\cal Q}$-Markov process) over an arbitrary $k$-tuple of sets in ${\cal C}$ and shows that,
if the transition system ${\cal Q}$ satisfies Assumption 1, then the definition will not depend on the extremal representations of these sets. 

\vspace{2mm}
{\bf Lemma 5.}
{\em Let ${\cal Q}:=(Q_{BB'})_{B \subseteq B'}$ be a transition
system satisfying Assumption 1. Let $(C_{1}, \ldots, C_{k})$ be a $k$-tuple of distinct
sets in ${\cal C}$ and suppose that each set $C_{i}$ admits two
extremal representations $C_{i}=A_{i} \verb2\2
\cup_{j=1}^{n_{i}}A_{ij}=A'_{i} \verb2\2
\cup_{j=1}^{m_{i}}A'_{ij}$. Let ${\cal A'}, {\cal A''}$ be the
minimal finite sub-semilattices of ${\cal A}$ which contain the
sets $A_{i},A_{ij}$, respectively $A'_{i},A'_{ij}$,
$\{B_{0}=\emptyset',B_{1}, \ldots, B_{n}\},
\{B'_{0}=\emptyset',B'_{1}, \ldots, B'_{m}\}$ two consistent
orderings of ${\cal A'}, {\cal A''}$ and $D_{j},D'_{l}$ the left
neighbourhoods of the sets $B_{j}, B'_{l}$ for $j=1,
\ldots,n;l=1, \ldots,m$. If each set $C_{i};i=1, \ldots,k$ can
be written as $C_{i}=\cup_{j \in J_{i}}D_{j}= \cup_{l \in
L_{i}}D'_{l}$ for some $J_{i} \subseteq \{1, \ldots,n\}, L_{i}
\subseteq \{1, \ldots,m\}$, then $$\int_{{\bf R}^{n+1}}
\prod_{i=1}^{k} I_{\Gamma_{i}}(\sum_{j \in
J_{i}}(x_{j}-x_{j-1})) Q_{\cup_{j=1}^{n-1}B_{j},
\cup_{j=1}^{n}B_{j}}(x_{n-1};dx_{n}) \ldots$$ $$Q_{B_{1},B_{1}
\cup B_{2}}(x_{1};dx_{2}) Q_{\emptyset'B_{1}}(x;dx_{1})
\mu(dx)=$$  $$\int_{{\bf R}^{m+1}} \prod_{i=1}^{k}
I_{\Gamma_{i}}(\sum_{l \in L_{i}}(y_{l}-y_{l-1}))
Q_{\cup_{l=1}^{m-1}B'_{l}, \cup_{l=1}^{m}B'_{l}}(y_{m-1};dy_{m})
\ldots$$ $$Q_{B'_{1},B'_{1} \cup
B'_{2}}(y_{1};dy_{2})Q_{\emptyset'B'_{1}}(y;dy_{1})\mu(dy)$$   
for every $\Gamma_{1}, \ldots, \Gamma_{k} \in {\cal B}({\bf
R})$,  with the convention $x_{0}=y_{0}=0$.  }   \vspace{2mm}

\noindent{\bf Proof}: Let $\tilde{\cal A}$ be the minimal
finite sub-semilattice of ${\cal A}$ determined by the sets in
${\cal A'}$ and ${\cal A''}$; clearly ${\cal A'} \subseteq
\tilde{\cal A}, {\cal A''} \subseteq \tilde{\cal A}$. Let ${\rm
ord}^{1}= \{E_{0}= \emptyset', E_{1}, \ldots, E_{N} \}$ and
${\rm ord}^{2}= \{E'_{0}= \emptyset', E'_{1}, \ldots, E'_{N} \}$
be two consistent orderings of $\tilde{\cal A}$ such that, if
$B_{j}=E_{i_{j}}; j=1, \ldots,n$ and $B'_{l}=E'_{k_{l}}; l=1,
\ldots,m$ for some indices $i_{1}<i_{2}< \ldots <i_{n}$,
respectively $k_{1}<k_{2}< \ldots <k_{m}$, then  
$$B_{1}=\cup_{p=1}^{i_{1}}E_{p}, \ B_{1} \cup B_{2}=
\cup_{p=1}^{i_{2}}E_{p}, \ldots, \cup_{j=1}^{n}B_{j}=
\cup_{p=1}^{i_{n}}E_{p}$$ 
$$B'_{1}= \cup_{q=1}^{k_{1}}E'_{q}, \ B'_{1} \cup B'_{2}=
\cup_{q=1}^{k_{2}}E'_{q}, \ldots, \cup_{l=1}^{m}B'_{l}=
\cup_{q=1}^{k_{m}}E'_{q}$$ 

Let $\pi$ be the permutation of $\{1,\ldots,N \}$ such that
$E_{p}=E'_{\pi(p)};p=1, \ldots,N$. Denote by $H_{p},H'_{q}$ the
left neighbourhoods of $E_{p},E'_{q}$ with respect to the
orderings ${\rm ord}^{1}, {\rm ord}^{2}$; clearly
$H_{p}=H'_{\pi(p)}, \ \forall p=1, \ldots, N$. 
Note that $D_{j}=(\cup_{v=1}^{j}B_{v}) \verb2\2
(\cup_{v=1}^{j-1}B_{v})=(\cup_{p=1}^{i_{j}}E_{p}) \verb2\2
(\cup_{p=1}^{i_{j-1}}E_{p})= \cup_{p=i_{j-1}+1}^{i_{j}}H_{p}; j=1,
\ldots,n$ and similarly $D'_{l}=
\cup_{q=k_{l-1}+1}^{k_{l}}H'_{q}; l=1, \ldots,m$. Hence
\begin{eqnarray*} C_{i} & = & \cup_{j \in J_{i}}D_{j} = \cup_{j
\in J_{i}} \cup_{p=i_{j-1}+1}^{i_{j}}H_{p} = \cup_{j \in J_{i}}
\cup_{p=i_{j-1}+1}^{i_{j}} H'_{\pi(p)} \\
      & = & \cup_{l \in L_{i}}D'_{l} = \cup_{l \in L_{i}}
\cup_{q=k_{l-1}+1}^{k_{l}}H'_{q}
\end{eqnarray*}
and we can conclude that
$$\{ \pi(p); p \in \cup_{j \in J_{i}} \{i_{j-1}+1,i_{j-1}+2,
\ldots, i_{j} \} \}= \cup_{l \in L_{i}} \{k_{l-1}+1,k_{l-1}+2,
\ldots, k_{l} \}$$
This implies that
$$\int_{{\bf R}^{N+1}} \prod_{i=1}^{k}I_{\Gamma_{i}}(\sum_{j \in
J_{i}} \sum_{p=i_{j-1}+1}^{i_{j}}(x_{p}-x_{p-1}))
Q_{\cup_{p=1}^{N-1}E_{p},\cup_{p=1}^{N}E_{p}}(x_{N-1};dx_{N}) 
\ldots $$
$$Q_{E_{1},E_{1} \cup E_{2}}(x_{1};dx_{2})
Q_{\emptyset' E_{1}}(x;dx_{1}) \mu(dx)$$
$$= \int_{{\bf R}^{N+1}} \prod_{i=1}^{k}I_{\Gamma_{i}}(\sum_{l
\in L_{i}} \sum_{q=k_{l-1}+1}^{k_{l}}(y_{q}-y_{q-1}))
Q_{\cup_{q=1}^{N-1}E'_{q},\cup_{q=1}^{N}E'_{q}}(y_{N-1};dy_{N}) 
\ldots $$
$$Q_{E'_{1},E'_{1} \cup E'_{2}}(y_{1};dy_{2})
Q_{\emptyset' E'_{1}}(y;dy_{1}) \mu(dy)$$
with the convention $x_{0}=y_{0}=0$. This gives us the desired
relationship, because $\sum_{p=i_{j-1}+1}^{i_{j}}
(x_{p}-x_{p-1})= x_{i_{j}}-x_{i_{j-1}}$, the left-hand side
integral collapses to an integral with respect to
$Q_{\cup_{j=1}^{n}B_{j}, \cup_{j=1}^{n-1}B_{j}}
(x_{i_{n-1}};dx_{i_{n}}) \ldots Q_{B_{1},B_{1} \cup
B_{2}}(x_{i_{1}}; dx_{i_{2}}) \linebreak Q_{\emptyset'
B_{1}}(x;dx_{i_{1}}) \mu(dx)$, and a similar phenomenon happens
in the right-hand side. $\Box$

We are now ready to prove the main result of this section, which
says that in fact the previous assumption is also sufficient
to construct a ${\cal Q}$-Markov process.

\vspace{2mm}
{\bf Theorem 1.}
{\em If $\mu$ is a probability measure on ${\bf R}$ and ${\cal
Q}:=(Q_{BB'})_{B \subseteq B'}$ is a transition system
which satisfies the consistency Assumption 1, then there exists
a ${\cal Q}$-Markov process $X:=(X_{A})_{A \in {\cal A}}$ with
initial distribution $\mu$. }  
\vspace{2mm}

\noindent{\bf Proof}: Let $({\bf R}^{\cal C}, {\cal R}^{\cal C}):= \prod_{C \in
{\cal C}} (R_{C}, {\cal R}_{C})$ where $(R_{C}, {\cal R}_{C}), C \in
{\cal C}$ represent ${\cal C}$ copies of the real space ${\bf R}$
with its Borel subsets. For each $k$-tuple $(C_{1}, \ldots,C_{k})$ of
distinct sets in ${\cal C}$ we will define a probability measure
$\mu_{C_{1} \ldots C_{k}}$ on $({\bf R}^{k}, {\cal R}^{k})$ such
that the system of all these probability measures satisfy the
following two consistency conditions:
\begin{description}
\item[(C1)] If $(C'_{1} \ldots C'_{k})$ is another ordering of 
the $k$-tuple $(C_{1} \ldots C_{k})$, say $C_{i}=C'_{\pi(i)}$, 
$\pi$ is a permutation of $\{1, \ldots,k \}$, then
$$\mu_{C_{1} \ldots C_{k}}(\Gamma_{1} \times \ldots \times
\Gamma_{k})=\mu_{C'_{1} \ldots C'_{k}}(\Gamma_{\pi^{-1}(1)} \times
\ldots \times \Gamma_{\pi^{-1}(k)})$$ for every $\Gamma_{1},
\ldots, \Gamma_{k} \in {\cal B}({\bf R})$. 

\item[(C2)] If $C_{1}, \ldots,C_{k},C_{k+1}$ are $k+1$ distinct
sets in ${\cal C}$, then
$$\mu_{C_{1} \ldots C_{k}}(\Gamma_{1} \times \ldots \times
\Gamma_{k})= \mu_{C_{1} \ldots C_{k} C_{k+1}}(\Gamma_{1} \times
\ldots \times \Gamma_{k} \times {\bf R})$$ 
for every $\Gamma_{1}, \ldots, \Gamma_{k} \in {\cal B}({\bf R})$.
\end{description}

By Kolmogorov's extension theorem there exists a probability measure
$P$ on $({\bf R}^{\cal C}, {\cal R}^{\cal C})$ such that the
coordinate-variable process $X:=(X_{C})_{C \in {\cal C}}$ defined by
$X_{C}(x):=x_{C}$ has the measures $\mu_{C_{1} \ldots C_{k}}$ as its
finite dimensional distributions. We will prove that the process $X$
has  an (almost surely) unique additive extension to ${\cal C}(u)$.
The
${\cal Q}$-Markov property of this process will follow from the way 
we choose its finite dimensional distributions $\mu_{C_{0}C_{1} 
\ldots C_{n}}$ over the left neighbourhoods of a finite
sub-semilattice, according to Proposition 5, (e). 

{\bf Step 1} (Construction of the finite dimensional distributions) 

\noindent We define $\mu_{\emptyset}:= \delta_{0}$. Let $(C_{1}, \ldots, C_{k})$ be a $k$-tuple of distinct nonempty sets
in ${\cal C}$ and $C_{i}=A_{i} \verb2\2 \cup_{j=1}^{n_{i}}A_{ij}; i=1,
\ldots,k$ some extremal representations. Let ${\cal
A'}$ be the minimal finite sub-semilattice of ${\cal A}$ which
contains the sets $A_{i},A_{ij}$, $\{B_{0}=\emptyset',B_{1},
\ldots, B_{n}\}$ a consistent ordering of ${\cal A'}$
and $D_{j}$ the left neighbourhood of the set $B_{j}$ for $j=1,
\ldots,n$. Define
$$\mu_{D_{0}D_{1} \ldots D_{n}} (\Gamma_{0} \times
\Gamma_{1} \times \ldots \times \Gamma_{n}):= \int_{{\bf R}^{n+1}}
I_{\Gamma_{0}}(x_{0}) I_{\Gamma_{1}}(x_{1}) 
\prod_{i=2}^{n} I_{\Gamma_{i}}(x_{i}-x_{i-1})$$
$$Q_{\cup_{j=1}^{n-1}B_{j}, \cup_{j=1}^{n}B_{j}}(x_{n-1};dx_{n})
\ldots Q_{B_{1},B_{1} \cup B_{2}}(x_{1};dx_{2})
Q_{\emptyset'B_{1}}(x_{0};dx_{1}) \mu(dx_{0})$$ 
for every $\Gamma_{0}, \Gamma_{1}, \ldots,\Gamma_{n} \in {\cal
B}({\bf R})$.

Say $C_{i}=\cup_{j \in J_{i}}D_{j}$ for some $J_{i}
\subseteq \{1, \ldots,n\};i=1, \ldots,k$, let $\alpha: {\bf R}^{n+1} \rightarrow
{\bf R}^{k}, \alpha(x_{0},x_{1}, \ldots,x_{n}):=(\sum_{j \in
J_{1}}x_{j}, \ldots,  \sum_{j \in J_{k}}x_{j})$ and define
$$\mu_{C_{1} \ldots C_{k}}:=  \mu_{D_{0}D_{1} \ldots D_{n}} \circ
\alpha^{-1}$$ 
The fact that $\mu_{C_{1} \ldots C_{k}}$ does not depend on the
ordering of the semilattice ${\cal A'}$ is a consequence of
Assumption 1.
The fact that $\mu_{C_{1} \ldots C_{k}}$ does not depend on
the extremal representations of the sets $C_{i}$ is also a
consequence of Assumption 1, using Lemma 5. Finally, we observe
that the finite dimensional distributions are additive.

{\bf Step 2} (Consistency condition (C1)) 

\noindent Let $(C'_{1} \ldots C'_{k})$ be another ordering of the
$k$-tuple $(C_{1} \ldots C_{k})$, with $C_{i}=C'_{\pi(i)}$, $\pi$
being a permutation of $\{1, \ldots,k \}$. 
Let $C_{i}=A_{i} \verb2\2 \cup_{j=1}^{n_{i}}A_{ij}; i=1,
\ldots,k$ be extremal representations, ${\cal
A'}$ the minimal finite sub-semilattice of ${\cal A}$ which
contains the sets $A_{i},A_{ij}$, $\{B_{0}=\emptyset',B_{1},
\ldots, B_{n}\}$ a consistent ordering of ${\cal A'}$
and $D_{j}$ the left neighbourhood of the set $B_{j}$ for each
$j=1, \ldots,n$. Say $C_{i}=\cup_{j \in J_{i}}D_{j},C'_{i}=
\cup_{j \in J'_{i}}D_{j}$ and let $\alpha(x_{0},x_{1},
\ldots,x_{n}):=(\sum_{j \in J_{1}}x_{j}, \ldots,  \sum_{j \in
J_{k}}x_{j}), \linebreak \beta(x_{0},x_{1}, \ldots,x_{n}):=(\sum_{j
\in J'_{1}}x_{j}, \ldots,  \sum_{j \in J'_{k}}x_{j})$. We have
$J_{i}=J'_{\pi(i)}$. Hence
\begin{eqnarray*}
 \mu_{C_{1} \ldots C_{k}}(\Gamma_{1} \times \ldots \times
\Gamma_{k}) &= &\mu_{D_{0}D_{1} \ldots D_{n}}(\alpha^{-1}(\Gamma_{1}
\times \ldots \times \Gamma_{k}))\\
&= &
 \mu_{D_{0}D_{1} \ldots D_{n}}(\beta^{-1}
(\Gamma_{\pi^{-1}(1)} \times \ldots \times \Gamma_{\pi^{-1}(k)}))\\
&=&
\mu_{C'_{1} \ldots C'_{k}}(\Gamma_{\pi^{-1}(1)} \times \ldots \times
\Gamma_{\pi^{-1}(k)}) 
\end{eqnarray*}
for every $\Gamma_{1},\ldots, \Gamma_{k} \in {\cal B}({\bf R})$.

{\bf Step 3} (Consistency Condition (C2))

\noindent Let $C_{1}, \ldots, C_{k}, C_{k+1}$ be $k+1$ distinct sets
in ${\cal C}$ and $C_{i}=A_{i} \verb2\2 \cup_{j=1}^{n_{i}}A_{ij}; i=1,
\ldots,k+1$ some extremal representations. Let ${\cal A'}, {\cal
A''}$ be the minimal finite sub-semilattices of ${\cal A}$ which
contain the sets $A_{i},A_{ij};i=1, \ldots,k;j=1, \ldots,n_{i}$,
respectively $A_{i},A_{ij};i=1, \ldots,k+1;j=1, \ldots,n_{i}$.
Clearly ${\cal A'} \subseteq {\cal A''}$. Using Comment 1.2 there
exists a consistent ordering $\{B_{0}= \emptyset',B_{1}, \ldots,B_{n}
\}$ of ${\cal A''}$ such that, if ${\cal A'}=\{B_{i_{0}}=
\emptyset',B_{i_{1}}, \ldots,B_{i_{m}} \}$ with $0=i_{0}< i_{1} <
\ldots < i_{m}$, then $\cup_{s=1}^{l}B_{i_{s}}=
\cup_{j=1}^{i_{l}}B_{j}$ for all $l=1, \ldots,m$. For each $j=1,
\ldots,n;l=1, \ldots,m$, let $D_{j},E_{l}$ be the left neighbourhoods
of $B_{j}$ in ${\cal A''}$, respectively of $B_{i_{l}}$ in ${\cal
A'}$ and note that  $E_{l}=\cup_{j=i_{l-1}+1}^{i_{l}}D_{j}$ for each
$l=1, \ldots,m$.  Let $\gamma(x_{0},x_{1}, \ldots,x_{n}):=(x_{0},
\sum_{j=1}^{i_{1}}x_{j}, \sum_{j=i_{1}+1}^{i_{2}}x_{j}, \ldots,
\sum_{j=i_{m-1}+1}^{i_{m}}x_{j})$ and for the moment suppose that 
\begin{equation}
\label{additive-type} \mu_{E_{0}E_{1} \ldots E_{m}}= \mu_{D_{0}D_{1}
\ldots D_{n}} \circ \gamma^{-1} 
\end{equation} 
Say $C_{i}=\cup_{j \in J_{i}}D_{j};i=1, \ldots,k+1$ for some $J_{i}
\subseteq \{1, \ldots, n \}$ and define $\alpha(x_{0},x_{1},
\linebreak \ldots,x_{n}):= (\sum_{j \in J_{1}}x_{j}, \ldots, \sum_{j
\in J_{k+1}}x_{j})$. On the other hand, if we say that  for each $i=1,
\ldots,k$ we have $C_{i}= \cup_{l \in L_{i}}E_{l}$ for some $L_{i}
\subseteq \{1, \ldots,m \}$, then $J_{i} = \cup_{l \in L_{i}}
\{i_{l-1}+1, i_{l-1}+2, \ldots,i_{l} \}$;
define $\beta(y_{0},y_{1}, \ldots,y_{m}):= (\sum_{l \in L_{1}}y_{l},
\linebreak \ldots, \sum_{l \in L_{k}}y_{l})$. Then   
\begin{eqnarray*}
 \mu_{C_{1}
\ldots C_{k}}(\Gamma_{1} \times \ldots \times \Gamma_{k})&=&
\mu_{E_{0}E_{1} \ldots E_{m}} (\beta^{-1}(\Gamma_{1} \times \ldots
\times \Gamma_{k}))\\
&= & \mu_{D_{0}D_{1} \ldots D_{n}} (\gamma^{-1}
(\beta^{-1}  (\Gamma_{1} \times \ldots \times \Gamma_{k}))) \\
&=&
\mu_{D_{0}D_{1} \ldots D_{n}} (\alpha^{-1}(\Gamma_{1} \times \ldots
\times \Gamma_{k} \times {\bf R}))\\
&= &
 \mu_{C_{1} \ldots C_{k}C_{k+1}}(\Gamma_{1}  \times \ldots \times
\Gamma_{k} \times {\bf R}) 
\end{eqnarray*} 
for every $\Gamma_{1}, \ldots,\Gamma_{k}
\in {\cal B}({\bf R})$.

In order to prove (\ref{additive-type}) let $\Gamma_{0},
\Gamma_{1}, \ldots, \Gamma_{m} \in {\cal B}({\bf R})$ be arbitrary.
Then $$\mu_{D_{0}D_{1} \ldots D_{n}}(\gamma^{-1}(\Gamma_{0} \times
\Gamma_{1} \times \ldots \times \Gamma_{m}))= $$
$$\int_{{\bf R}^{n+1}}
I_{\Gamma_{0} \times \Gamma_{1}
\times \ldots \times \Gamma_{m}}(\gamma(x_{0},x_{1},x_{2}-x_{1},
\ldots,x_{n}-x_{n-1}))$$
$$Q_{\cup_{j=1}^{n-1}B_{j},\cup_{j=1}^{n}B_{j}}(x_{n-1};dx_{n}) \ldots
Q_{B_{1},B_{1} \cup
B_{2}}(x_{1};dx_{2})Q_{\emptyset'B_{1}}(x_{0};dx_{1}) \mu(dx_{0})$$
Note that $\gamma(x_{0},x_{1},x_{2}-x_{1}, \ldots,x_{n}-x_{n-1})=
(x_{0},x_{i_{1}},x_{i_{2}}-x_{i_{1}}, \ldots,x_{i_{m}}-x_{i_{m-1}})$
and hence the integrand does not depend on the variables $x_{j}, j
\not \in \{i_{0},i_{1}, \ldots,i_{m} \}$. By the definition of the
transition system, the above integral collapses to  $$\int_{{\bf
R}^{m+1}} I_{\Gamma_{0}}(x_{0}) I_{\Gamma_{1}}(x_{i_{1}})
\prod_{l=2}^{m}I_{\Gamma_{l}}(x_{i_{l}}-x_{i_{l-1}})
Q_{\cup_{j=1}^{i_{m-1}}B_{j},\cup_{j=1}^{i_{m}}B_{j}}(x_{i_{m-1}};
dx_{i_{m}}) \ldots$$ $$Q_{\cup_{j=1}^{i_{1}}B_{j},
\cup_{j=1}^{i_{2}}B_{j}}(x_{i_{1}};dx_{i_{2}}) Q_{\emptyset',
\cup_{j=1}^{i_{1}}B_{j}}(x_{0};dx_{i_{1}}) \mu(dx_{0})$$ which is
exactly the definition of $\mu_{E_{0}E_{1} \ldots E_{m}}(\Gamma_{0}
\times \Gamma_{1} \times \ldots \times \Gamma_{m})$ because
$\cup_{j=1}^{i_{l}}B_{j}= \cup_{s=1}^{l}B_{i_{s}}$ for every $l=1,
\ldots,m$. This concludes the proof of (\ref{additive-type}).

{\bf Step 4} (Almost Sure Additivity of the Canonical Process)

\noindent We will show that the canonical process $X$ on the
space $({\bf R}^{\cal C}, {\cal R}^{\cal C})$ has an (almost surely)
unique  additive extension to ${\cal C}(u)$ (with respect to the
probability measure $P$ given by Kolmogorov's extension theorem).
Let $C,C_{1},\ldots,C_{k} \in {\cal C}$ be such that
$C=\cup_{i=1}^{k}C_{i}$, and suppose that $C_{i}=A_{i} \verb2\2
\cup_{j=1}^{n_{i}}A_{ij}; i=1, \ldots,k$ are extremal
representations. Let ${\cal A'}$ be the minimal finite 
sub-semilattice of ${\cal A}$ which contains the sets $A_{i},A_{ij}$,
$\{B_{0}=\emptyset', B_{1}, \ldots, B_{n}\}$ a consistent ordering of
${\cal A'}$ and $D_{j}$ the left neighbourhood of $B_{j}$. Assume
that each $C_{i}= \cup_{j \in J_{i}}D_{j}$ for some $J_{i} \subseteq
\{1, \ldots, n \}$. Because the finite dimensional distributions of
$X$ were chosen in an additive way, we have $X_{C_{i}}= \sum_{j \in
J_{i}}X_{D_{j}}$ a.s., $X_{C_{i_{1}} \cap C_{i_{2}}}= \sum_{j \in
I_{i_{1}} \cap I_{i_{2}}}X_{D_{j}}$ a.s., $\ldots$ ,
$X_{\cap_{i=1}^{k}C_{i}}=\sum_{j \in \cap_{i=1}^{k}I_{i}}X_{D_{j}}$
a.s., $X_{C}= \sum_{j \in \cup_{i=1}^{k} J_{i}}X_{D_{j}}$ a.s. Hence
$\sum_{i=1}^{k}X_{C_{i}}- \sum_{i_{1}<i_{2}}X_{C_{i_{1}} \cap
C_{i_{2}}}+ \ldots +(-1)^{k}X_{\cap_{i=1}^{k}C_{i}}= \sum_{j \in
\cup_{i=1}^{k} J_{i}}X_{D_{j}}=X_{C}$ a.s. 

$\Box$

Finally we will translate the preceding result in terms of
flows. Let $\mu$ be an arbitrary probability measure on ${\bf R}$.

For each finite sub-semilattice ${\cal A'}$ and for each consistent
ordering ord= $\{A_{0}= \emptyset',A_{1}, \ldots,A_{n} \}$ of ${\cal
A'}$ pick one simple flow $f:=f_{{\cal A'}, {\rm ord}}$ which
connects the sets of ${\cal A'}$ in the sense of the ordering
ord. Let ${\cal S}$ be the collection of all the simple flows
$f_{{\cal A'}, {\rm ord}}$ and $\{ {\cal Q}^{f}:=(Q_{st}^{f})_{s<t};
f \in {\cal S}\}$ a collection of one-dimensional transition systems
indexed by ${\cal S}$.  

The next assumption will provide a necessary and sufficient
condition which will allow us to reconstruct a set-indexed
transition system ${\cal Q}$ from a class of one-dimensional
transition systems $\{Q^f\}$. It requires that whenever we have
two simple flows $f,g \in {\cal S}$ such that
$f(s)=g(u),f(t)=g(v)$ for some $s<t,u<v$, we must have
$Q_{st}^{f}=Q_{uv}^{g}$. 

\vspace{2mm}
{\bf Assumption 2.}
{\em If ord1$=\{A_{0}=\emptyset',A_{1}, \ldots,A_{n}\}$ and 
ord2$=\{A_{0}=\emptyset',A'_{1}, \linebreak \ldots, A'_{m}\}$ are two
consistent orderings of some finite semilattices ${\cal A'}, {\cal
A''}$ such that $\cup_{j=1}^{n}A_{j}= \cup_{j=1}^{m}A'_{j}$,
$\cup_{j=1}^{k}A_{j}= \cup_{j=1}^{l}A'_{j}$ for some $k<n,l<m$, and
we denote $f:=f_{{\cal A'},ord1}$, $g:=f_{{\cal A''},ord2}$ with
$f(t_{i})= \cup_{j=1}^{i}A_{j}$, $g(u_{i})= \cup_{j=1}^{i}A'_{j}$,
then }
$$Q_{0t_{1}}^{f}=Q_{0u_{1}}^{g}
\ {\rm and} \ Q_{t_{k}t_{n}}^{f}=Q_{u_{l}u_{m}}^{g}$$

The following assumption is easily seen to be equivalent to
Assumption 1.

\vspace{2mm}
{\bf Assumption 3.}
{\em If ord$1=\{A_{0}=\emptyset',A_{1}, \ldots,A_{n}\}$ and 
ord$2=\{A_{0}=\emptyset',A'_{1}, \linebreak \ldots, A'_{n}\}$ are two
consistent orderings of the same finite semilattice ${\cal A'}$ with 
$A_{i}=A'_{\pi(i)}$, where $\pi$ is a permutation of 
$\{1, \ldots,n\}$ with $\pi(1)=1$, and we denote $f:=f_{{\cal
A'},ord1}$, $g:=f_{{\cal A'},ord2}$ with $f(t_{i})=
\cup_{j=1}^{i}A_{j}$, $g(u_{i})= \cup_{j=1}^{i}A'_{j}$, then 
\begin{equation}
\label{integr-cond-1}   
\int_{{\bf R}^{n+1}} I_{\Gamma_{0}}(x_{0})
I_{\Gamma_{1}}(x_{1})\prod_{i=2}^{n} I_{\Gamma_{i}}(x_{i}-x_{i-1})
Q_{t_{n-1}t_{n}}^{f}(x_{n-1};dx_{n}) \ldots 
\end{equation}
$$Q_{t_{1}t_{2}}^{f}(x_{1};dx_{2}) Q_{0t_{1}}^{f}(x_{0};dx_{1})
\mu(dx_{0})=$$   
$$ \int_{{\bf R}^{n+1}}
I_{\Gamma_{0}}(y_{0})I_{\Gamma_{1}}(y_{1})
\prod_{i=2}^{n}I_{\Gamma_{i}}(y_{\pi(i)}-y_{\pi(i)-1})Q_{u_{n-1}u_{n}}^{g}
(y_{n-1};dy_{n}) \ldots$$
$$ Q_{u_{1}u_{2}}^{g}(y_{1};dy_{2})
Q_{0u_{1}}^{g}(y_{0};dy_{1}) \mu(dy_{0})$$  
for every $\Gamma_{0}, \Gamma_{1}, \ldots, \Gamma_{n} \in {\cal
B}({\bf R})$. }   
\vspace{2mm} 

The following result is immediate. 

\vspace{2mm}
{\bf Corollary 1.}
{\em If $\mu$ is a probability measure on ${\bf R}$ and $\{ {\cal
Q}^{f}:=(Q_{st}^{f})_{s<t}; \linebreak f \in {\cal S}\}$ is a
collection of one-dimensional transition systems which satisfies
the matching Assumption 2 and the consistency Assumption 3,
then there exist a set-indexed transition
system ${\cal Q}:=(Q_{BB'})_{B \subseteq B'}$ and a ${\cal
Q}$-Markov process $X:=(X_{A})_{A \in {\cal A}}$ with initial
distribution $\mu$, such that $\forall f \in {\cal S}$, 
$X^{f}:=(X_{f(t)})_{t}$ is ${\cal Q}^{f}$-Markov. }
\vspace{2mm}

\noindent{\bf Proof}: Let $B,B' \in {\cal A}(u)$ be such that $B \subseteq
B'$. Let $B=\cup_{j=1}^{m}A'_{j}, B'= \cup_{l=1}^{p}A''_{l}$ be some
extremal representations, ${\cal A'}$ the minimal finite
sub-semilattice of ${\cal A}$ which contains the sets $A'_{j},
A''_{l}$ and ord=$\{A_{0}= \emptyset',A_{1}, \ldots,A_{n}\}$ a
consistent ordering of ${\cal A'}$ such that $B=
\cup_{j=1}^{k}A_{j}$ and $B'= \cup_{j=1}^{n}A_{j}$. 
Denote $f:=f_{{\cal A'}, {\rm ord}}$ with $f(t_{i})=
\cup_{j=1}^{i}A_{j}$. We define $Q_{BB'}:=Q_{t_{k}t_{n}}^{f}$. 

The definition of $Q_{BB'}$ does not depend on the extremal
representations of $B,B'$, because of Assumption 2.

Finally it is not hard to see that the family ${\cal Q}:=(Q_{BB'})_{B
\subseteq B'}$ is a transition system which satisfies 
Assumption 1. $\Box$

\section{The Generator}\label{Generator}

In this section we will make use of flows to introduce the generator
of a ${\cal Q}$-Markov process. Corollary 1 allows us
to characterize a ${\cal Q}$-Markov process by a class $\{{\cal
Q}^{f}; f \in {\cal S} \}$ of one-dimensional transition systems.
These transition sytems in turn are characterized by their
corresponding generators. This observation permits us to define a
generator for a ${\cal Q}$-Markov process as a class of generators
indexed by a suitable class of flows. The generator completely characterizes the distribution 
of a ${\cal Q}$-Markov process. We will also determine necessary and sufficient conditions
that will ensure that a family of one-dimensional generators, indexed by a collection of simple flows, is the generator of a ${\cal Q}$-Markov process.

Let $B({\bf R})$ be the Banach space of all bounded measurable
functions $h: {\bf R} \rightarrow {\bf R}$ with the supremum norm.

Let ${\cal Q}:=(Q_{BB'})_{B \subseteq B'}$ be a
transition system and $X:=(X_{A})_{A \in {\cal A}}$ a
${\cal Q}$-Markov process with initial distribution $\mu$. For each
finite sub-semilattice ${\cal A'}$ and for each consistent ordering
ord= $\{A_{0}= \emptyset',A_{1}, \ldots,A_{n} \}$ of ${\cal A'}$ pick
one simple flow $f:=f_{{\cal A'}, {\rm ord}}$ which connects the sets
of ${\cal A'}$ in the sense of the ordering ord. Let ${\cal S}$ be the
collection of all the simple flows $f_{{\cal A'}, {\rm ord}}$. 

For each $f \in {\cal S}$, let ${\cal T}^{f}:=(T_{st}^{f})_{s<t}$ be
the semigroup associated to the transition system ${\cal Q}^{f}$ and
${\cal G}_{s}^{f},{\cal G}_{s}^{*f}$ the backward, respectively the
forward generator of the process $X^{f}$ at time $s$, with
domains ${\cal D}({\cal G}_{s}^{f}), {\cal D}({\cal G}_{s}^{*f})$.
We will make the usual assumption that all the domains ${\cal
D}({\cal G}_{s}^{f})$ and ${\cal D}({\cal G}_{s}^{*f})$ have a
common subspace ${\cal D}$, which is dense in $B({\bf R})$, such
that for every $s<t$, $T_{st}^{f}({\cal D}) \subseteq {\cal D}$
and for every $h \in {\cal D}$, the function $T_{st}^{f}h$ is
strongly continuously differentiable with respect to $s$ and
$t$ with $$- ({{\partial^{-}} \over {\partial
r}}T_{rs}^{f}h)|_{r=s}= ({{\partial^{+}} \over {\partial
t}}T_{st}^{f}h)|_{t=s} \ \ \forall s$$ 
The operator ${\cal G}_{s}^{f}={\cal G}_{s}^{*f}$ defined on
${\cal D}$, will be called {\em the generator} of the process
$X^{f}$ at time $s$. A consequence of Kolmogorov-Feller equations is
that 
\begin{equation}
\label{integral-expression}
T_{st}^{f}h-h= \int_{s}^{t}{\cal G}_{v}^{f}T_{vt}^{f} h \ dv, \
\forall h \in {\cal D}
\end{equation}

\vspace{2mm}
{\bf Definition 5.}
{\em The collection $\{ {\cal G}^{f}:=({\cal G}_{s}^{f})_{s}; f \in
{\cal S} \}$, where ${\cal G}_{s}^{f}$ is the generator of the
one-dimensional process $X^{f}$ at time $s$, is called the {\bf
generator} of the set-indexed process $X$. }
\vspace{2mm}

{\bf Corollary 2.}
{\em The generator of a ${\cal Q}$-Markov process determines its distribution.}
\vspace{2mm}

\noindent {\bf Proof}: Since ${\cal G}^{f}$ determines ${\cal T}^{f}$ and ${\cal Q}^{f}$, the result is an immediate consequence of Corollary 1. $\Box$

\vspace{2mm}
{\bf Examples 4.}
\begin{enumerate}
\item The generator of a process with independent increments, which is
stochastically continuous and weakly differentiable on every simple
flow $f \in {\cal S}$, is given by  
\begin{eqnarray*}  
\lefteqn{({\cal G}_{s}^{f} h)(x)=(\gamma_{s}^{f})' h'(x) + {{1} \over
{2}} (\Lambda_{s}^{f})' h''(x) + } \\
 & & \int_{\{|y|>1\}}(h(x+y)-h(x))
(\Pi_{s}^{f})'(dx) + \\
 & & \int_{\{|y| \leq 1\}} (h(x+y)-h(x)-y h'(x))
(\Pi_{s}^{f})'(dx) 
\end{eqnarray*} 
where $\gamma_{t}^{f},\Lambda_{t}^{f},\Pi_{t}^{f}(dy)$ are
respectively, the translation function, the variance measure,
and the L\'{e}vy measure of the process on the flow $f$ and
$(\gamma_{s}^{f})', (\Lambda_{s}^{f})',(\Pi_{s}^{f})'(dy)$ denote the
derivatives at $s$ of these functions. The domain of ${\cal
G}_{s}^{f}$ contains the dense subspace $C_{b}^{2}$ of twice
continuously differentiable functions $h: {\bf R} \rightarrow {\bf
R}$ such that $h,h',h''$ are bounded.

\item The generator of the empirical process of size $n$,
corresponding to a probability measure $F$ which has the
property that $\forall f \in {\cal S} \ F \circ f$ is differentiable,
is given by      
$$({\cal G}_{s}^{f}h) \left({{k} \over {n}} \right)=
(n-k) \left( h \left({{k+1} \over {n}} \right) -h \left({{k} \over {n}}
\right) \right) {{(F \circ f)'(s)} \over {1-F \circ f(s)}}, \ k<n,$$
where $(F \circ f)'(s)$ denotes the derivative at $s$ of $F \circ f$.
The domain of ${\cal G}_{s}^{f}$ is the space of all finite arrays
$h:=(h(k))_{k=0, \ldots,n}$.

\item The generator of the Dirichlet process with parameter measure $\alpha$ which has the property that $\forall f \in {\cal S} \ \alpha \circ f$ is differentiable, is given by
$$({\cal G}_{s}^{f}h)(x)=(\alpha \circ f)'(s) \int_{0}^{1-x}{{h(x+y)-h(x)} \over {y}} \left( {{1-x-y} \over {1-x}} \right)^{\alpha(f(s)^{c})-1}dy$$
where $(\alpha \circ f)'(s)$ denotes the derivative at $s$ of $\alpha \circ f$.
The domain of ${\cal G}_{s}^{f}$ contains the dense subspace $C_{b}^{1}$ of  
continuously differentiable functions $h: {\bf R} \rightarrow {\bf
R}$ such that $h,h'$ are bounded.

\end{enumerate}
\vspace{2mm}

The goal is to find the conditions that have to be satisfied by
the collection $\{ {\cal G}^{f}:=({\cal G}_{s}^{f})_{s}; f \in 
{\cal S} \}$ of generators such that the associated collection $\{
{\cal Q}^{f}; f \in {\cal S} \}$ of one-parameter transition systems
satisfies the matching Assumption 2 and the consistency
Assumption 3. By invoking Corollary 1, we will be able to
conclude next that there exists a set-indexed transition system
${\cal Q}$ and a ${\cal Q}$-Markov process $X:=(X_{A})_{A \in
{\cal A}}$ with initial distribution $\mu$, whose generator is
exactly the collection $\{{\cal G}^{f}:=({\cal G}_{s}^{f})_{s};
f \in {\cal S} \}$.

Using the integral equation (\ref{integral-expression}), we will first
give the equivalent form in terms of generators of the matching
Assumption 2.

\vspace{2mm}
{\bf Assumption 4.}
{\em If ord$1=\{A_{0}=\emptyset',A_{1}, \ldots,A_{n}\}$ and 
ord$2=\{A_{0}=\emptyset',A'_{1}, \linebreak \ldots, A'_{m}\}$ are
two consistent orderings of some finite semilattices ${\cal
A'}, {\cal A''}$,  such that $\cup_{j=1}^{n}A_{j}=
\cup_{j=1}^{m}A'_{j}, \cup_{j=1}^{k}A_{j}=
\cup_{j=1}^{l}A'_{j}$ for some $k<n,l<m$, and we denote
$f:=f_{{\cal A'}, ord1}$, $g:=f_{{\cal A''},ord2}$ with
$f(t_{i})= \cup_{j=1}^{i}A_{j}$, $g(u_{i})=
\cup_{j=1}^{i}A'_{j}$, then for every $h \in {\cal D}$
$$\int_{0}^{t_{1}} {\cal G}_{v}^{f}T_{vt_{1}}^{f}h dv
=\int_{0}^{u_{1}} {\cal G}_{v}^{g}T_{vu_{1}}^{g}h dv \ {\rm and} \
\int_{t_{k}}^{t_{n}} {\cal G}_{v}^{f}T_{vt_{n}}^{f}h dv
=\int_{u_{l}}^{u_{m}} {\cal G}_{v}^{g}T_{vu_{m}}^{g}h dv$$ }

We will need the following notational convention. 

\noindent If $Q_{1}(x_{1};dx_{2}), Q_{2}(x_{2};dx_{3}), \ldots,
Q_{n-1}(x_{n-1};dx_{n})$ are transition probabilities on ${\bf
R}$, $T_{1}, T_{2}, \ldots, T_{n-1}$ are their associated bounded
linear operators, and $h: {\bf R}^{n} \rightarrow {\bf R}$ is a
bounded measurable function, then     $$T_{1}T_{2} \ldots T_{n-1} \
[h(x_{1},x_{2}, \ldots,x_{n})](x_{1}) : \stackrel{\rm def}{=}$$ 
$$\int_{{\bf R}^{n-1}} h(x_{1},x_{2}, \ldots,x_{n})
Q_{n-1}(x_{n-1};dx_{n}) \ldots Q_{2}(x_{2};dx_{3})
Q_{1}(x_{1};dx_{2})$$
If in addition, ${\cal G}$ is a
linear operator on $B({\bf R})$ with domain ${\cal D}({\cal
G})$, and the function $h$ is chosen such that for every
$x_{1}, \ldots,x_{k} \in {\bf R}$, the function
$$\int_{{\bf R}^{n-k}} h(x_{1}, x_{2}, \ldots,x_{n})
Q_{n-1}(x_{n-1};dx_{n}) \ldots
Q_{k+1}(x_{k+1};dx_{k+2})Q_{k}(\cdot; dx_{k+1})$$ is
in ${\cal D}({\cal G})$, then    
$$T_{1}T_{2} \ldots T_{k-1} {\cal G} T_{k}
\ldots T_{n-1} \ [h(x_{1},x_{2}, \ldots,x_{n})]
(x_{1}) : \stackrel{\rm def}{=}$$ $$\int_{{\bf R}^{k-1}} ({\cal
G} \int_{{\bf R}^{n-k}} h(x_{1}, x_{2}, \ldots,x_{n})
Q_{n-1}(x_{n-1};dx_{n}) \ldots Q_{k+1}(x_{k+1};dx_{k+2})$$
$$Q_{k}(\cdot;dx_{k+1}))(x_{k})Q_{k-1}(x_{k-1};dx_{k}) \ldots
Q_{2}(x_{2};dx_{3})Q_{1}(x_{1};dx_{2}).$$
    
In order to understand the consistency Assumption 3 and to see how 
we can express it in terms of the generators, we consider a simple example.
Let ${\cal A'}$ be a finite semilattice consisting of $6$ sets,  ord$1=\{A_{0}=\emptyset',A_{1}, \ldots,A_{6}\}$ and 
ord$2=\{A_{0}=\emptyset',A'_{1}, \ldots, A'_{6}\}$ two
consistent orderings of ${\cal A'}$ with
$A_{i}=A'_{\pi(i)}$, where $\pi$ is the permutation $(1)(24)(36)(5)$. Denote with $C_{i},C'_{i}$ the left neighbourhoods of $A_{i}$ respectively $A'_{i}$.

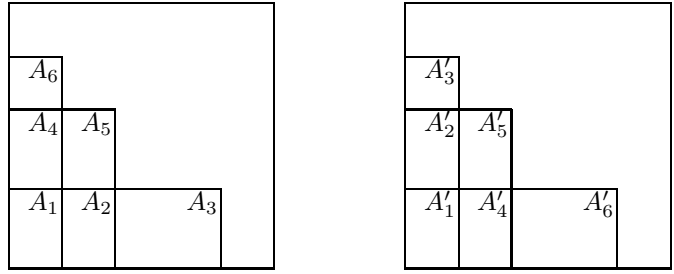
\begin{center}
\begin{picture}(250,100)
   \put(0,0){\framebox(100,100)}
      \put(20,0){\line(0,1){80}}
      \put(40,0){\line(0,1){60}}
      \put(80,0){\line(0,1){30}}
      \put(0,30){\line(1,0){80}}
      \put(0,60){\line(1,0){40}}
      \put(0,80){\line(1,0){20}}
      \put(7,22){$A_{1}$}
      \put(27,22){$A_{2}$}
      \put(67,22){$A_{3}$}
      \put(7,52){$A_{4}$}
      \put(27,52){$A_{5}$}
      \put(7,72){$A_{6}$}

\put(150,0){\framebox(100,100)}
     \put(170,0){\line(0,1){80}}
      \put(190,0){\line(0,1){60}}
      \put(230,0){\line(0,1){30}}
      \put(150,30){\line(1,0){80}}
      \put(150,60){\line(1,0){40}}
      \put(150,80){\line(1,0){20}}
      \put(157,22){$A'_{1}$}
      \put(177,22){$A'_{4}$}
      \put(217,22){$A'_{6}$}
      \put(157,52){$A'_{2}$}
      \put(177,52){$A'_{5}$}
      \put(157,72){$A'_{3}$}
\end{picture}
\end{center}

Set $f:=f_{{\cal A'},ord1}$, $g:=f_{{\cal A'},ord2}$ with $f(t_{i})=
\cup_{j=1}^{i}A_{j}$, $g(u_{i})= \cup_{j=1}^{i}A'_{j}$ for $i=1, \ldots,6$.
Suppose that a ${\cal Q}$-Markov process $(X_{A})_{A \in {\cal A}}$ exists; denote $X:=X^{f}$ and $Y:=X^{g}$. We have
$$X_{t_{2}}-X_{t_{1}}=Y_{u_{4}}-Y_{u_{3}}; \ \ X_{t_{3}}-X_{t_{2}}=Y_{u_{6}}-Y_{u_{5}}; \ \
X_{t_{4}}-X_{t_{3}}=Y_{u_{2}}-Y_{u_{1}}$$
$$X_{t_{5}}-X_{t_{4}}=Y_{u_{5}}-Y_{u_{4}}; \ \ X_{t_{6}}-X_{t_{5}}=Y_{u_{3}}-Y_{u_{2}}$$

Using these equalities we will first find the relationship between $T_{t_{i-1}t_{i}}^{f}$ and $T_{u_{\pi(i)-1},u_{\pi(i)}}^{f}$ and then between the generators $({\cal G}_{w}^{f})_{w \in [t_{i-1},t_{i}]}$ and $({\cal G}_{v}^{g})_{v \in [u_{\pi(i)-1},u_{\pi(i)}]}$ for $i=2, \ldots,6$.

For $i=2$ we have $P[X_{t_{2}}-X_{t_{1}} \in \Gamma|X_{t_{1}}=x_{1}] =
P[Y_{u_{4}}-Y_{u_{3}} \in \Gamma|Y_{u_{1}}=x_{1}]$
which leads us to the equation: $\forall h \in B({\bf R})$
$$\int_{\bf R} h(x_{2}) Q_{t_{1}t_{2}}^{f}(x_{1};dx_{2})=
\int_{{\bf R}^{2}} h(x_{1}+y_{4}-y_{3}) Q_{u_{3} u_{4}}^{g}(y_{3};dy_{4}) Q_{u_{1}
u_{3}}^{g}(x_{1};dy_{3})$$

\noindent Using the above notational convention, this can be written as
\begin{equation}
\label{special-permut:i=2}
(T_{t_{1}t_{2}}^{f}h)(x_{1})=T_{u_{1}u_{3}}^{g}T_{u_{3}u_{4}}^{g}[h(x_{1}+y_{4}-y_{3})](x_{1}).
\end{equation}

\noindent Let $h \in {\cal D}$ and subtract $h(x_{1})$ from both sides of this equation. Using (\ref{integral-expression}) and Fubini's theorem, we get 
\begin{equation}
\label{special-permut-gen:i=2}
\int_{t_{1}}^{t_{2}}({\cal G}_{w}^{f}T_{wt_{2}}^{f}h)(x_{1})dw=\int_{u_{3}}^{u_{4}}T_{u_{1}u_{3}}^{g}{\cal G}_{v}^{g}T_{vu_{4}}^{g}[h(x_{1}+y_{4}-y_{3})](x_{1})dv
\end{equation}
Note that (\ref{special-permut:i=2}) and (\ref{special-permut-gen:i=2}) are equivalent (since ${\cal D}$ is dense in $B({\bf R})$).

\vspace{3mm}

For $i=3$, $P[X_{t_{2}}-X_{t_{1}} \in \Gamma_{1}, X_{t_{3}}-X_{t_{2}} \in \Gamma_{2}|X_{t_{1}}=x_{1}] =
P[Y_{u_{4}}-Y_{u_{3}} \in \Gamma_{1},Y_{u_{6}}-Y_{u_{5}} \in \Gamma_{2} |Y_{u_{1}}=x_{1}]$ which leads us to the equation: $\forall h_{1},h_{2} \in B({\bf R})$
\begin{equation}
\label{special-permut:i=3}
T_{u_{1}u_{3}}^{g}T_{u_{3}u_{4}}^{g}[h_{1}(y_{4}-y_{3})(T_{t_{2}t_{3}}^{f}h_{2})(x_{1}+y_{4}-y_{3})](x_{1})
\end{equation}
$$=T_{u_{1}u_{3}}^{g}T_{u_{3}u_{4}}^{g}T_{u_{4}u_{5}}^{g}T_{u_{5}u_{6}}^{g}[h_{1}(y_{4}-y_{3})h_{2}(x_{1}+y_{4}-y_{3}+y_{6}-y_{5})](x_{1})$$

\noindent Let $h_{1},h_{2} \in {\cal D}$ and subtract $T_{u_{1}u_{3}}^{g}T_{u_{3}u_{4}}^{g}[h_{1}(y_{4}-y_{3})h_{2}(x_{1}+y_{4}-y_{3})](x_{1})$ from both sides of this equation. On the left-hand 
side we get 
$$\int_{t_{2}}^{t_{3}} T_{u_{1}u_{3}}^{g}T_{u_{3}u_{4}}^{g}[h_{1}(y_{4}-y_{3})
({\cal G}_{w}^{f}T_{wt_{3}}^{f}h_{2})(x_{1}+y_{4}-y_{3})](x_{1})$$

\noindent On the right-hand side we have
$$T_{u_{1}u_{3}}^{g}T_{u_{3}u_{4}}^{g}T_{u_{4}u_{5}}^{g} [(T_{u_{5}u_{6}}^{g}h'(x_{1},y_{3},y_{4},y_{5}, \cdot))(y_{5})-h(x_{1},y_{3},y_{4})](x_{1})$$
with
$h'(x_{1},y_{3},y_{4},y_{5},y_{6}):=h_{1}(y_{4}-y_{3})h_{2}(x_{1}+y_{4}-y_{3}+y_{6}-y_{5})$ and
$h(x_{1},y_{3},y_{4})=h_{1}(y_{4}-y_{3})h_{2}(x_{1}+y_{4}-y_{3})$, which becomes
$$\int_{u_{5}}^{u_{6}}T_{u_{1}u_{3}}^{g}T_{u_{3}u_{4}}^{g}T_{u_{4}u_{5}}^{g} {\cal G}_{v}^{g}T_{vu_{6}}^{g}[h_{1}(y_{4}-y_{3})h_{2}(x_{1}+y_{4}-y_{3}+y_{6}-y_{5})](x_{1})$$
because $h'(x_{1},y_{3},y_{4},y_{5},y_{5})=h(x_{1},y_{3},y_{4})$.

\noindent Therefore, the equivalent form of (\ref{special-permut:i=3}) in terms of the generators is
\begin{equation}
\label{special-permut-gen:i=3}
\int_{t_{2}}^{t_{3}} T_{u_{1}u_{3}}^{g}T_{u_{3}u_{4}}^{g}[h_{1}(y_{4}-y_{3})
({\cal G}_{w}^{f}T_{wt_{3}}^{f}h_{2})(x_{1}+y_{4}-y_{3})](x_{1})
\end{equation}
$$=\int_{u_{5}}^{u_{6}}T_{u_{1}u_{3}}^{g}T_{u_{3}u_{4}}^{g}T_{u_{4}u_{5}}^{g} {\cal G}_{v}^{g}T_{vu_{6}}^{g}[h_{1}(y_{4}-y_{3})h_{2}(x_{1}+y_{4}-y_{3}+y_{6}-y_{5})](x_{1})$$

Continuing in the same manner for $i=4$ we will get the necessary relationship between 
$T_{t_{3}t_{4}}^{f}$ and $T_{u_{1}u_{2}}^{g}$. (When we write down this relationship we have to specify the ordering of $\pi(2)-1, \pi(2), \pi(3)-1, \pi(3), \pi(4)-1, \pi(4)$.) This relationship will have an equivalent form in terms of the generators $({\cal G}_{w}^{f})_{w \in [t_{3},t_{4}]}$ and $({\cal G}_{v})_{v \in [u_{1},u_{2}]}$.

At the end of this procedure we discover a collection of $5$ necessary relationships that have to be satisfied by the generators ${\cal G}^{f}$ and ${\cal G}^{g}$ of the process. A very important fact is that these relationships are also sufficient, i.e. if they hold then the finite dimensional distribution of the process over the semilattice ${\cal A'}$ is `invariant' under the two orderings {\em ord1} and {\em ord2}.

\vspace{3mm}

We return now to the general context. Let $\mu_{t_{1}}^{f}$ be the probability measure defined by $\mu_{t_{1}}^{f}( \Gamma):= \int_{\bf R}Q_{0t_{1}}^{f}(x; \Gamma)
\mu(dx)$. The following result is crucial and its non-trivial proof can be found
in the appendix.

\vspace{2mm}
{\bf Lemma 6.}
{\em Let ord$1=\{A_{0}=\emptyset',A_{1}, \ldots,A_{n}\}$ and 
ord$2=\{A_{0}=\emptyset',A'_{1}, \linebreak \ldots, A'_{n}\}$ be two
consistent orderings of the same finite semilattice ${\cal A'}$ with
$A_{i}=A'_{\pi(i)}$, where $\pi$ is a permutation of $\{1, \ldots,n\}$
with $\pi(1)=1$, and set $f:=f_{{\cal
A'},ord1}$, $g:=f_{{\cal A'},ord2}$ with $f(t_{i})=
\cup_{j=1}^{i}A_{j}$, $g(u_{i})= \cup_{j=1}^{i}A'_{j}$.

Suppose that $Q_{0t_{1}}^{f}=Q_{0u_{1}}^{g}$. The following
statements are equivalent:
\begin{description}
\item[(a)] the integral condition (\ref{integr-cond-1}) of
the consistency Assumption 3 holds;

\item[(b)] for each $i=2, \ldots,n$, if we denote with 
$l_{1} \leq l_{2} \leq \ldots \leq l_{2(i-1)}$ the increasing
ordering of the values $\pi(2)-1,\pi(2), \pi(3)-1,\pi(3), \ldots,
\pi(i)-1,\pi(i)$, and with $p$ the index for which
$\pi(i)-1=l_{p-1},\pi(i)=l_{p}$, then for every $h_{2},
\ldots,h_{i} \in {\cal D}$ and for $\mu_{t_{1}}^{f}$-almost
all $x_{1}$ 

{\bf (b1)} if $p=2(i-1)$ we have
$$\int_{t_{i-1}}^{t_{i}} T_{u_{1}u_{l_{1}}}^{g}
T_{u_{l_{1}}u_{l_{2}}}^{g} \ldots
T_{u_{l_{p-3}}u_{l_{p-2}}}^{g}$$
$$[\prod_{j=2}^{i-1}h_{j}(y_{\pi(j)}
-y_{\pi(j)-1})
({\cal G}_{w}^{f}T_{wt_{i}}^{f} h_{i})
(x_{1}+ \sum_{j=2}^{i-1} (y_{\pi(j)}-y_{\pi(j)-1}))](x_{1})dw$$
$$=\int_{u_{l_{p-1}}}^{u_{l_{p}}} 
T_{u_{1}u_{l_{1}}}^{g} T_{u_{l_{1}}u_{l_{2}}}^{g} \ldots
T_{u_{l_{p-3}}u_{l_{p-2}}}^{g} T_{u_{l_{p-2}}u_{l_{p-1}}}^{g}
{\cal G}_{v}^{g} T_{vu_{l_{p}}}^{g}$$
$$[\prod_{j=2}^{i-1}h_{j}(y_{\pi(j)}-y_{\pi(j)-1})h_{i}(x_{1}+
\sum_{j=2}^{i}(y_{\pi(j)}-y_{\pi(j)-1}))](x_{1})dv;$$

{\bf (b2)} if $p< 2(i-1)$ we have
$$\int_{t_{i-1}}^{t_{i}} 
T_{u_{1}u_{l_{1}}}^{g}  \ldots
T_{u_{l_{p-3}}u_{l_{p-2}}}^{g}
T_{u_{l_{p-2}}u_{l_{p+1}}}^{g}T_{u_{l_{p+1}}u_{l_{p+2}}}^{g}
\ldots T_{u_{l_{2(i-1)-1}}u_{l_{2(i-1)}}}^{g}$$
$$[\prod_{j=2}^{i-1}h_{j}(y_{\pi(j)}-y_{\pi(j)-1})
({\cal G}_{w}^{f}T_{wt_{i}}^{f} h_{i})
(x_{1}+\sum_{j=2}^{i-1} (y_{\pi(j)}-y_{\pi(j)-1}))](x_{1})dw$$
$$= \int_{u_{l_{p-1}}}^{u_{l_{p}}} 
T_{u_{1}u_{l_{1}}}^{g} T_{u_{l_{1}}u_{l_{2}}}^{g} \ldots
T_{u_{l_{p-3}}u_{l_{p-2}}}^{g}T_{u_{l_{p-2}}u_{l_{p-1}}}^{g}
{\cal G}_{v}^{g}T_{vu_{l_{p}}}^{g}$$
$$T_{u_{l_{p}}u_{l_{p+1}}}^{g} T_{u_{l_{p+1}}u_{l_{p+2}}}^{g}
\ldots T_{u_{l_{2(i-1)-1}}u_{l_{2(i-1)}}}^{g} \{
\prod_{j=2}^{i-1}h_{j}(y_{\pi(j)}-y_{\pi(j)-1}) \cdot$$  
$$[h_{i}(x_{1}+ \sum_{j=2}^{i}(y_{\pi(j)}-y_{\pi(j)-1}))-
h_{i}(x_{1}+ \sum_{j=2}^{i-1}(y_{\pi(j)}-y_{\pi(j)-1})) ]
\}(x_{1}) dv.$$ 
\end{description} }

The following assumption is the equivalent form in terms of
generators of the consistency Assumption 3.

\vspace{2mm}
{\bf Assumption 5.}
{\em If ord$1=\{A_{0}=\emptyset',A_{1}, \ldots,A_{n}\}$ and 
ord$2=\{A_{0}=\emptyset',A'_{1}, \linebreak \ldots, A'_{n}\}$ are two
consistent orderings of the same finite semilattice ${\cal A'}$ with 
$A_{i}=A'_{\pi(i)}$, where $\pi$ is a permutation of 
$\{1, \ldots,n\}$ with $\pi(1)=1$, and we denote $f:=f_{{\cal
A'},ord1}$, $g:=f_{{\cal A'},ord2}$ with $f(t_{i})=
\cup_{j=1}^{i}A_{j}$, $g(u_{i})= \cup_{j=1}^{i}A'_{j}$, then the
generators ${\cal G}^{f}$ and ${\cal G}^{g}$ satisfy condition (b)
stated in Lemma 6. } 
\vspace{2mm}

Clearly Assumptions 4 and 5 are necessary conditions satisfied by  
the generator of any set-indexed ${\cal Q}$-Markov process. The following result is an immediate consequence of Corollary 1 which says that in fact, these two assumptions are also sufficient
for the construction of the process.

\vspace{2mm}
{\bf Theorem 2.}
{\em Let $\mu$ be a probability measure on ${\bf R}$ and
$\{ {\cal G}^{f}:=({\cal G}_{s}^{f})_{s}; f \in {\cal S} \}$ a
collection of families of linear operators on $B({\bf R})$ such 
that each operator ${\cal G}_{s}^{f}$ is defined on a dense
subspace ${\cal D}$ of $B({\bf R})$ and is the generator at
time $s$ of a semigroup ${\cal T}^{f}:=(T_{st}^{f})_{s<t}$
associated with a transition system ${\cal
Q}^{f}:=(Q_{st}^{f})_{s<t}$. If the family $\{ {\cal G}^{f}; f \in
{\cal S} \}$ satisfies the matching Assumption 4 and the
consistency Assumption 5, then there exist a set-indexed
transition system ${\cal Q}:=(Q_{BB'})_{B,B' \in {\cal A}(u); B
\subseteq B'}$ and a ${\cal Q}$-Markov process $X:=(X_{A})_{A
\in {\cal A}}$ with initial distribution $\mu$, whose generator
is exactly the collection $\{ {\cal G}_{s}^{f}; f \in {\cal S}
\}$. }

\vspace{3mm}

\footnotesize{{\em Acknowledgement.} The authors greatly appreciate
the thoughtful comments and suggestions of Professor Donald. A. Dawson
regarding the use of the integral equation (\ref{integral-expression}) for proving the equivalence 
of statements (a) and (b) in Lemma 6.}

\normalsize{

\appendix

\section{Proof of Lemma 6}

We will prove the desired equivalence by means of an intermediate condition {\bf (b')}. We will show that {\bf (a) $\Leftrightarrow$ (b')} and {\bf (b') $\Leftrightarrow$ (b)}. Here is this condition.

\vspace{3mm}

{\bf (b')} {\em For each $i=2, \ldots,n$, if we denote with $l_{1} \leq
l_{2} \leq \ldots \leq l_{2(i-1)}$ the increasing ordering of the
values $\pi(2)-1,\pi(2), \pi(3)-1,\pi(3), \ldots, \pi(i)-1,\pi(i)$, and
with $p$ the index for which $\pi(i)-1=l_{p-1},\pi(i)=l_{p}$, then for
every $h_{2}, \ldots,h_{i} \in B({\bf R})$ and for
$\mu_{t_{1}}^{f}$-almost all $x_{1}$}
\begin{equation}
\label{equi}
T_{u_{1}u_{l_{1}}}^{g}
T_{u_{l_{1}}u_{l_{2}}}^{g} \ldots T_{u_{l_{p-3}}u_{l_{p-2}}}^{g}
T_{u_{l_{p-2}}u_{l_{p+1}}}^{g}T_{u_{l_{p+1}}u_{l_{p+2}}}^{g}
\ldots T_{u_{l_{2(i-1)-1}}u_{l_{2(i-1)}}}^{g}
\end{equation}
$$[\prod_{j=2}^{i-1}h_{j}(y_{\pi(j)}-y_{\pi(j)-1})
(T_{t_{i-1}t_{i}}^{f}h_{i})(x_{1}+\sum_{j=2}^{i-1}
(y_{\pi(j)}-y_{\pi(j)-1}))](x_{1})$$  
$$= T_{u_{1}u_{l_{1}}}^{g}
T_{u_{l_{1}}u_{l_{2}}}^{g} \ldots
T_{u_{l_{2(i-1)-1}}u_{l_{2(i-1)}}}^{g}$$
$$[\prod_{j=2}^{i-1}h_{j}(y_{\pi(j)}-y_{\pi(j)-1})
h_{i}(x_{1}+\sum_{j=2}^{i}(y_{\pi(j)}-y_{\pi(j)-1}))](x_{1});$$

{\bf Proof of (a) $\Rightarrow$ (b')}: Let $X$ be a ${\cal
Q}^{f}$-Markov process and $Y$ a ${\cal Q}^{g}$-Markov process
with the same initial distribution $\mu$. Since
$Q_{0t_{1}}^{f}=Q_{0u_{1}}^{g}$, both $X_{t_{1}}$ and
$Y_{u_{1}}$ have the same distribution $\mu_{t_{1}}^{f}$. The
integral condition (\ref{integr-cond-1}) of the consistency
Assumption 6 is equivalent to saying that for 
$\mu_{t_{1}}^{f}$-almost all $x_{1}$, the conditional
distribution of  $(X_{t_{2}}-X_{t_{1}},
\ldots,X_{t_{n}}-X_{t_{n-1}})$ given $X_{t_{1}}=x_{1}$ coincide
with the conditional distribution of 
$(Y_{u_{\pi(2)}}-Y_{u_{\pi(2)-1}},\ldots,Y_{u_{\pi(n)}}-Y_{u_{\pi(n)-1}})$
given $Y_{u_{1}}=x_{1}$.

For $i=2$ we will use the following relationship:
for every $\Gamma_{2} \in {\cal B}({\bf R})$ and for
$\mu_{t_{1}}^{f}$-almost all $x_{1}$ 
$$P[X_{t_{2}}-X_{t_{1}} \in \Gamma_{2}|X_{t_{1}}=x_{1}] =
P[Y_{u_{\pi(2)}}-Y_{u_{\pi(2)-1}} \in \Gamma_{2}|
Y_{u_{1}}=x_{1}]$$
Using the Markov property, the left-hand side is $\int_{\bf R}
I_{\Gamma_{2}+x_{1}}(x_{2})Q_{t_{1} t_{2}}^{f}(x_{1};dx_{2})$, whereas
on the right-hand side we have $$\int_{{\bf R}^{2}} 
I_{\Gamma_{2}}(y_{u_{\pi(2)}}-y_{u_{\pi(2)-1}}) Q_{u_{\pi(2)-1}
u_{\pi(2)}}^{g}(y_{\pi(2)-1}; dy_{\pi(2)}) Q_{u_{1} u_{\pi(2)-1}}^{g}
(x_{1}; dy_{\pi(2)-1})$$   
By a monotone class argument we can conclude that for every $h_{2} \in
B({\bf R})$ and for $\mu_{t_{1}}^{f}$-almost all $x_{1}$
\begin{equation}
\label{equ2}
\int_{\bf R} h_{2}(x_{2}) Q_{t_{1}t_{2}}^{f}(x_{1};dx_{2})=
\int_{{\bf R}^{2}} h_{2}(x_{1}+y_{u_{\pi(2)}}-y_{u_{\pi(2)-1}}) 
\end{equation}
$$Q_{u_{\pi(2)-1} u_{\pi(2)}}^{g}(y_{\pi(2)-1};dy_{\pi(2)}) Q_{u_{1}
u_{\pi(2)-1}}^{g}(x_{1};dy_{\pi(2)-1})$$
which is the desired relationship for $i=2$.

For $i=3$ we will use the following relationship: for every
$\Gamma_{2}, \Gamma_{3} \in {\cal B}({\bf R})$ and for
$\mu_{t_{1}}^{f}$-almost all $x_{1}$ \begin{eqnarray*}  
\lefteqn{P [X_{t_{2}}-X_{t_{1}} \in
\Gamma_{2}, X_{t_{3}}-X_{t_{2}} \in \Gamma_{3}|X_{t_{1}}=x_{1}]=
} \\
 & & P[ Y_{u_{\pi(2)}}-Y_{u_{\pi(2)-1}} \in
\Gamma_{2},Y_{u_{\pi(3)}}-Y_{u_{\pi(3)-1}} \in
\Gamma_{3}|Y_{u_{1}}=x_{1}] 
\end{eqnarray*}

The left-hand side can be written as
$$\int_{{\bf R}^{2}} 
I_{\Gamma_{2} +x_{1}}(x_{2})I_{\Gamma_{3}+x_{2}}(x_{3}) Q_{t_{2}
t_{3}}^{f}(x_{2};dx_{3}) Q_{t_{1} t_{2}}^{f} (x_{1};dx_{2})$$
which becomes
$$\int_{{\bf R}^{3}} 
I_{\Gamma_{2}}(y_{\pi(2)}-y_{\pi(2)-1})
I_{\Gamma_{3}+ x_{1}+y_{\pi(2)}-y_{\pi(2)-1}}(x_{3})
Q_{t_{2} t_{3}}^{f}(x_{1}+y_{\pi(2)}-y_{\pi(2)-1};dx_{3})$$
$$Q_{u_{\pi(2)-1} u_{\pi(2)}}^{g}(y_{\pi(2)-1};
dy_{\pi(2)})^{g} Q_{u_{1} u_{\pi(2)-1}}^{g}
(x_{1};dy_{\pi(2)-1})$$
using equation (\ref{equ2}).

The right-hand side can be written as
$$\int_{{\bf R}^{4}} 
I_{\Gamma_{2}}(y_{\pi(2)}-y_{\pi(2)-1})
I_{\Gamma_{3}}(y_{\pi(3)}-y_{\pi(3)-1})Q_{u_{l_{3}},u_{l_{4}}}^{g}
(y_{l_{3}};dy_{l_{4}})
Q_{u_{l_{2}},u_{l_{3}}}^{g} (y_{l_{2}};dy_{l_{3}})$$
$$Q_{u_{l_{1}},u_{l_{2}}}^{g} (y_{l_{1}};dy_{l_{2}})
Q_{u_{1},u_{l_{1}}}^{g} (x_{1};dy_{l_{1}})$$ 
where $l_{1} \leq l_{2} \leq l_{3} \leq l_{4}$  is the increasing
ordering of the values $\pi(2)-1,\pi(2), \pi(3)
\linebreak -1,\pi(3)$.

By a monotone class argument we can conclude that for every
$h_{2},h_{3} \in B({\bf R})$ and for  $\mu_{t_{1}}^{f}$-almost all
$x_{1}$
\begin{equation}
\label{equ3}
\int_{{\bf R}^{3}} h_{2}(y_{\pi(2)}-y_{\pi(2)-1}) h_{3}(x_{3})
Q_{t_{2} t_{3}}^{f}(x_{1}+y_{\pi(2)}-y_{\pi(2)-1};dx_{3})
\end{equation}
$$Q_{u_{\pi(2)-1} u_{\pi(2)}}^{g}(y_{\pi(2)-1};dy_{\pi(2)})
Q_{u_{1} u_{\pi(2)-1}}^{g}(x_{1};dy_{\pi(2)-1})$$
$$= \int_{{\bf R}^{4}} h_{2}(y_{\pi(2)}-y_{\pi(2)-1}) h_{3}(x_{1}+
y_{\pi(2)}-y_{\pi(2)-1}+y_{\pi(3)}-y_{\pi(3)-1})
Q_{u_{l_{3}},u_{l_{4}}}^{g}(y_{l_{3}};dy_{l_{4}})$$
$$Q_{u_{l_{2}},u_{l_{3}}}^{g}(y_{l_{2}};dy_{l_{3}}) 
Q_{u_{l_{1}},u_{l_{2}}}^{g} (y_{l_{1}};dy_{l_{2}})
Q_{u_{1},u_{l_{1}}}^{g}(x_{1};dy_{l_{1}})$$
which is the desired relationship for $i=3$.

The inductive argument will be omitted since it is identical,
but notationally complex. 

{\bf Proof of (b')} $\Rightarrow$ {\bf (a)}: Let $\Gamma_{0},\Gamma_{1}, \ldots,
\Gamma_{n} \in {\cal B}({\bf R})$ be arbitrary. Using the fact
that $Q_{0t_{1}}^{f}=Q_{0u_{1}}^{g}$ and equation
(\ref{equ2}),  the left-hand side of the integral
condition (\ref{integr-cond-1}) becomes   $$\int_{{\bf R}^{n+2}}
I_{\Gamma_{0}}(y_{0}) I_{\Gamma_{1}}(y_{1}) I_{\Gamma_{2}}(y_{\pi(2)}-
y_{\pi(2)-1}) I_{\Gamma_{3}+y_{1}+y_{\pi(2)}-y_{\pi(2)-1}}(x_{3})
I_{\Gamma_{4}+x_{3}}(x_{4}) \ldots$$ $$I_{\Gamma_{n}+x_{n-1}}(x_{n})
Q_{t_{n-1}t_{n}}^{f}(x_{n-1};dx_{n}) \ldots
Q_{t_{3}t_{4}}^{f}(x_{3};dx_{4})
Q_{t_{2}t_{3}}^{f}(y_{1}+y_{\pi(2)}-y_{\pi(2)-1}; dx_{3})$$
$$Q_{u_{\pi(2)-1} u_{\pi(2)}}^{g}(y_{\pi(2)-1};
dy_{\pi(2)})Q_{u_{1} u_{\pi(2)-1}}^{g} (y_{1};dy_{\pi(2)-1})
Q_{0u_{1}}^{g}(y_{0};dy_{1}) \mu(dy_{0})$$ 
which in turn can be written as   
$$\int_{{\bf R}^{n+3}} I_{\Gamma_{0}}(y_{0}) I_{\Gamma_{1}}(y_{1})
I_{\Gamma_{2}}(y_{\pi(2)}- y_{\pi(2)-1})
I_{\Gamma_{3}} (y_{\pi(3)}-y_{\pi(3)-1})$$
$$I_{\Gamma_{4}+ y_{1}+ \sum_{j=2}^{3}
(y_{\pi(j)}-y_{\pi(j)-1})}(x_{4}) \ldots
I_{\Gamma_{n}+x_{n-1}}(x_{n})
Q_{t_{n-1}t_{n}}^{f}(x_{n-1};dx_{n}) \ldots$$
$$Q_{t_{3}t_{4}}^{f}( y_{1}+ \sum_{j=2}^{3}
(y_{\pi(j)}-y_{\pi(j)-1});dx_{4})
Q_{u_{l_{3}},u_{l_{4}}}^{g} (y_{l_{3}};dy_{l_{4}})
Q_{u_{l_{2}},u_{l_{3}}}^{g} (y_{l_{2}};dy_{l_{3}})$$
$$Q_{u_{l_{1}},u_{l_{2}}}^{g} (y_{l_{1}};dy_{l_{2}})
Q_{u_{1},u_{l_{1}}}^{g} (y_{1};dy_{l_{1}})
Q_{0u_{1}}^{g}(y_{0};dy_{1}) \mu(dy_{0})$$  using equation
(\ref{equ3}), where $l_{1} \leq l_{2} \leq l_{3} \leq l_{4}$ 
is the increasing ordering of the values $\pi(2)-1,\pi(2),
\pi(3)-1,\pi(3)$.

Continuing in the same manner at the last step we will get
exactly the desired right-hand side of equation
(\ref{integr-cond-1}), since the increasing ordering of the values
$\pi(2)-1, \pi(2), \ldots, \pi(n)-1,\pi(n)$ is exactly $1 \leq 2
\leq \ldots \leq n$. 

{\bf Proof of (b') $\Leftrightarrow$ (b)}: The basic ingredient will be
equation (\ref{integral-expression}), which gives the integral
expression of a semigroup in terms of its generator.

Since ${\cal D}$ is dense we can assume that the functions $h_{2},
\ldots,h_{i}$ are in ${\cal D}$ in the expression given by {\bf (b')}.
Subtract $$T_{u_{1}u_{l_{1}}}^{g}
T_{u_{l_{1}}u_{l_{2}}}^{g} \ldots T_{u_{l_{p-3}}u_{l_{p-2}}}^{g}
T_{u_{l_{p-2}}u_{l_{p+1}}}^{g}T_{u_{l_{p+1}}u_{l_{p+2}}}^{g}
\ldots T_{u_{l_{2(i-1)-1}}u_{l_{2(i-1)}}}^{g}$$
$$[\prod_{j=2}^{i-1}h_{j}(y_{\pi(j)}-y_{\pi(j)-1})
h_{i}(x_{1}+\sum_{j=2}^{i-1}(y_{\pi(j)}-y_{\pi(j)-1}))](x_{1})$$
from both sides of this expression.

On the left-hand side we will have 
$$T_{u_{1}u_{l_{1}}}^{g}
T_{u_{l_{1}}u_{l_{2}}}^{g} \ldots T_{u_{l_{p-3}}u_{l_{p-2}}}^{g}
T_{u_{l_{p-2}}u_{l_{p+1}}}^{g}T_{u_{l_{p+1}}u_{l_{p+2}}}^{g}
\ldots T_{u_{l_{2(i-1)-1}}u_{l_{2(i-1)}}}^{g}$$
$$[\prod_{j=2}^{i-1}h_{j}(y_{\pi(j)}-y_{\pi(j)-1})
(T_{t_{i-1}t_{i}}^{f}h_{i}-h_{i})(x_{1}+\sum_{j=2}^{i-1}
(y_{\pi(j)}-y_{\pi(j)-1}))](x_{1})$$
which can be written as 
$$\int_{t_{i-1}}^{t_{i}}
T_{u_{1}u_{l_{1}}}^{g} T_{u_{l_{1}}u_{l_{2}}}^{g} \ldots
T_{u_{l_{p-3}}u_{l_{p-2}}}^{g}
T_{u_{l_{p-2}}u_{l_{p+1}}}^{g}T_{u_{l_{p+1}}u_{l_{p+2}}}^{g}
\ldots T_{u_{l_{2(i-1)-1}}u_{l_{2(i-1)}}}^{g}$$ 
$$[\prod_{j=2}^{i-1}h_{j}(y_{\pi(j)}-y_{\pi(j)-1})
({\cal G}_{w}^{f}T_{wt_{i}}^{f} h_{i})
(x_{1}+\sum_{j=2}^{i-1} (y_{\pi(j)}-y_{\pi(j)-1}))](x_{1})dw$$ 
using the integral expression (\ref{integral-expression}) and
Fubini's theorem.

On the right-hand side we have
$$T_{u_{1}u_{l_{1}}}^{g}
T_{u_{l_{1}}u_{l_{2}}}^{g} \ldots T_{u_{l_{p-3}}u_{l_{p-2}}}^{g}
T_{u_{l_{p-2}}u_{l_{p-1}}}^{g}$$
$$\{ T_{u_{l_{p-1}}u_{l_{p}}}^{g} T_{u_{l_{p}}u_{l_{p+1}}}^{g}
T_{u_{l_{p+1}}u_{l_{p+2}}}^{g} \ldots
T_{u_{l_{2(i-1)-1}}u_{l_{2(i-1)}}}^{g}$$
$$[\prod_{j=2}^{i-1}h_{j}(y_{\pi(j)}-y_{\pi(j)-1})
h_{i}(x_{1}+
\sum_{j=2}^{i}(y_{\pi(j)}-y_{\pi(j)-1}))](y_{l_{p-1}})-$$
$$-T_{u_{l_{p-1}}u_{l_{p+1}}}^{g} T_{u_{l_{p+1}}u_{l_{p+2}}}^{g}
\ldots T_{u_{l_{2(i-1)-1}}u_{l_{2(i-1)}}}^{g}$$
$$[\prod_{j=2}^{i-1}h_{j}(y_{\pi(j)}-y_{\pi(j)-1})
h_{i}(x_{1}+
\sum_{j=2}^{i-1}(y_{\pi(j)}-y_{\pi(j)-1}))](y_{l_{p-1}})
\}(x_{1})$$

If we denote
$$h'(x_{1},y_{l_{1}},
\ldots,y_{l_{p-1}},y_{l_{p}}):=$$ $$\int_{{\bf R}^{2(i-1)-p}}
\prod_{j=2}^{i-1}h_{j}(y_{\pi(j)}-y_{\pi(j)-1})
h_{i}(x_{1}+ \sum_{j=2}^{i}(y_{\pi(j)}-y_{\pi(j)-1}))$$
$$Q_{u_{l_{2(i-1)-1}}u_{l_{2(i-1)}}}^{g}(y_{l_{2(i-1)-1}};
dy_{l_{2(i-1)}})\ldots Q_{u_{l_{p+1}}u_{l_{p+2}}}^{g}
(y_{l_{p+1}};dy_{l_{p+2}})$$
$$Q_{u_{l_{p}}u_{l_{p+1}}}^{g}(y_{l_{p}};dy_{l_{p+1}})$$ and
$$h(x_{1},y_{l_{1}}, \ldots,y_{l_{p-1}}):=$$
$$\int_{{\bf R}^{2(i-1)-p}} 
\prod_{j=2}^{i-1}h_{j}(y_{\pi(j)}-y_{\pi(j)-1})
h_{i}(x_{1}+ \sum_{j=2}^{i-1}(y_{\pi(j)}-y_{\pi(j)-1}))$$
$$Q_{u_{l_{2(i-1)-1}}u_{l_{2(i-1)}}}^{g}(y_{l_{2(i-1)-1}};
dy_{l_{2(i-1)}})\ldots Q_{u_{l_{p+1}}u_{l_{p+2}}}^{g}
(y_{l_{p+1}};dy_{l_{p+2}})$$
$$Q_{u_{l_{p-1}}u_{l_{p+1}}}^{g}(y_{l_{p-1}};dy_{l_{p+1}})$$
then the right-hand side becomes
$$T_{u_{1}u_{l_{1}}}^{g}
T_{u_{l_{1}}u_{l_{2}}}^{g} \ldots T_{u_{l_{p-3}}u_{l_{p-2}}}^{g}
T_{u_{l_{p-2}}u_{l_{p-1}}}^{g}$$
$$[(T_{u_{l_{p-1}}u_{l_{p}}}^{g}h'(x_{1},y_{l_{1}},\ldots,
y_{l_{p-1}},\cdot))(y_{l_{p-1}})-h(x_{1},y_{l_{1}},
\ldots, y_{l_{p-1}})](x_{1})$$ 
$$=T_{u_{1}u_{l_{1}}}^{g} T_{u_{l_{1}}u_{l_{2}}}^{g} \ldots
T_{u_{l_{p-3}}u_{l_{p-2}}}^{g} T_{u_{l_{p-2}}u_{l_{p-1}}}^{g}$$
$$[h'(x_{1},y_{l_{1}}, \ldots, y_{l_{p-1}},
y_{l_{p-1}})-h(x_{1},y_{l_{1}}, \ldots, 
y_{l_{p-1}})+$$
$$+ \int_{u_{l_{p-1}}}^{u_{l_{p}}} ({\cal G}_{v}^{g}
T_{vu_{l_{p}}}^{g}h'(x_{1},y_{l_{1}}, \ldots, 
y_{l_{p-1}},\cdot))(y_{l_{p-1}})dv](x_{1})$$
since $h'(x_{1},y_{l_{1}}, \ldots, y_{l_{p-1}}, 
\cdot) \in {\cal D}$ for every $x_{1},y_{l_{1}},
\ldots, y_{l_{p-1}}$. We have two cases:

\vspace{3mm}

{\bf Case 1)} If $p=2(i-1)$, then all the integrals with respect
to $Q_{u_{l_{p}}u_{l_{p+1}}}^{g},
\linebreak Q_{u_{l_{p+1}}u_{l_{p+2}}}^{g},
\ldots,Q_{u_{l_{2(i-1)-1}}u_{l_{2(i-1)}}}^{g}$ disappear in the
preceding expressions. Hence $h'(x_{1},y_{l_{1}},
\ldots,y_{l_{p-1}},y_{l_{p-1}})=
h(x_{1},y_{l_{1}}, \ldots, y_{l_{p-1}})$ and the
result follows.

\vspace{3mm}

{\bf Case 2)} If $p< 2(i-1)$, then
$$h'(x_{1},y_{l_{1}}, \ldots,y_{l_{p-1}},
y_{l_{p-1}})=(T_{u_{l_{p}}u_{l_{p+1}}}^{g}
H(x_{1},y_{l_{1}}, \ldots, y_{l_{p-1}},
\cdot))(y_{l_{p-1}})$$
$$h(x_{1},y_{l_{1}}, \ldots,y_{l_{p-1}})=
(T_{u_{l_{p-1}}u_{l_{p+1}}}^{g} H(x_{1},y_{l_{1}},
 \ldots,y_{l_{p-1}}, \cdot))(y_{l_{p-1}})$$ where
$$H(x_{1},y_{l_{1}}, \ldots, y_{l_{p-1}},
 y_{l_{p+1}}):=$$
$$\int_{{\bf R}^{2(i-1)-p-1}} 
\prod_{j=2}^{i-1}h_{j}(y_{\pi(j)}-y_{\pi(j)-1})
h_{i}(x_{1}+ \sum_{j=2}^{i-1}(y_{\pi(j)}-y_{\pi(j)-1}))$$
$$Q_{u_{l_{2(i-1)-1}}u_{l_{2(i-1)}}}^{g}(y_{l_{2(i-1)-1}};
dy_{l_{2(i-1)}}) \ldots Q_{u_{l_{p+1}}u_{l_{p+2}}}^{g}
(y_{l_{p+1}};dy_{l_{p+2}})$$

To simplify the notation we will omit the arguments of the
function $H$. Hence
$$h'(x_{1},y_{l_{1}}, \ldots, y_{l_{p-1}},
y_{l_{p-1}})- h(x_{1},y_{l_{1}}, \ldots, 
y_{l_{p-1}})$$
$$=(T_{u_{l_{p}}u_{l_{p+1}}}^{g} H)(y_{l_{p-1}})-
[T_{u_{l_{p-1}}u_{l_{p}}}^{g}(T_{u_{l_{p}}u_{l_{p+1}}}^{g}H)]
(y_{l_{p-1}})$$
$$= - \int_{u_{l_{p-1}}}^{u_{l_{p}}} ({\cal
G}_{v}^{g}T_{vu_{l_{p}}}^{g}
(T_{u_{l_{p}}u_{l_{p+1}}}^{g}H))(y_{l_{p-1}})dv$$
since $H(x_{1},y_{l_{1}}, \ldots, y_{l_{p-1}}, 
\cdot) \in {\cal D}$, for every $x_{1},y_{l_{1}},
\ldots, y_{l_{p-1}}$.

Using Fubini's theorem, the right-hand side becomes 
$$\int_{u_{l_{p-1}}}^{u_{l_{p}}}  
T_{u_{1}u_{l_{1}}}^{g} T_{u_{l_{1}}u_{l_{2}}}^{g} \ldots
T_{u_{l_{p-3}}u_{l_{p-2}}}^{g} T_{u_{l_{p-2}}u_{l_{p-1}}}^{g}$$
$$[({\cal G}_{v}^{g}T_{vu_{l_{p}}}^{g}
U(x_{1},y_{l_{1}}, \ldots,y_{l_{p-1}},
\cdot))(y_{l_{p-1}})](x_{1})dv$$
where
$$U(x_{1},y_{l_{1}}, \ldots,y_{l_{p-1}},y_{l_{p}}):=$$
$$=h'(x_{1},y_{l_{1}}, \ldots,y_{l_{p-1}},
y_{l_{p}})-(T_{u_{l_{p}}u_{l_{p+1}}}^{g}H(x_{1},
y_{l_{1}},y_{l_{2}},\ldots, y_{l_{p-1}}, \cdot))(y_{l_{p}})$$
$$= \int_{{\bf R}^{2(i-1)-p}} 
\prod_{j=2}^{i-1}h_{j}(y_{\pi(j)}-y_{\pi(j)-1})$$  
$$[h_{i}(x_{1}+ \sum_{j=2}^{i}(y_{\pi(j)}-y_{\pi(j)-1}))-
h_{i}(x_{1}+ \sum_{j=2}^{i-1}(y_{\pi(j)}-y_{\pi(j)-1}))]$$
$$Q_{u_{l_{2(i-1)-1}}u_{l_{2(i-1)}}}^{g}(y_{l_{2(i-1)-1}};
dy_{l_{2(i-1)}}) \ldots Q_{u_{l_{p}}u_{l_{p+1}}}^{g}
(y_{l_{p}};dy_{l_{p+1}})$$
This concludes the proof. $\Box$

 }

\end{document}